\documentclass[10pt,a4paper]{article}%
\usepackage[centertags]{amsmath}
\usepackage{amsfonts}
\usepackage{amssymb}
\usepackage{amsthm}
\usepackage{epsfig}
\usepackage{setspace}
\usepackage{ae}
\usepackage{eucal}
\usepackage[usenames]{color}%
\usepackage{graphicx}

\theoremstyle{plain}
\newtheorem{thm}{Theorem}[section]
\newtheorem{lem}[thm]{Lemma}

\theoremstyle{definition}

\theoremstyle{remark}
\newtheorem{exmp}[thm]{Example}
\newtheorem{rem}[thm]{Remark}

\renewcommand{\div}{\operatorname{div}}

\begin{document}

\title{Some interesting class of integrable partial differential equation systems (related to models of fluid mechanics)}
\author{J\"org Kampen }
\maketitle

\begin{abstract}
We determine a considerable class of nonlinear partial differential equation systems, related to a large class of fluid mechanical models, which have global regular solutions. Uniqueness is not a direct general  consequence of this method.  The solution scheme uses local-time representations in terms of probability densities and applies in strong function spaces. We observe that strong Gaussian damping of scaled models can offset possible growth of subhomogeneous nonlinear terms, where the spatial scaling parameter $r>1$  indicates the deviation from a strong semigroup property in strong norms. The application to the Navier Stokes equation is considered as an example of the scheme in detail. Although this construction of global solution branches works for a large class of models (and can be extended to models with variable viscosity terms and even to models with highly degenerate diffusions), we remark that a direct application to randomized models (for example average models or models with white noise) is not possible in general, as energy can be transported to higher frequencies caused by stochastic mixing. We also compare with results obtained from singularity analysis and CKN theory, and consider generalisations to models with variable diffusions or highly degenerate diffusions.   
\end{abstract}


2010 Mathematics Subject Classification.  35D35, 35Q35
\section{Definition of a integrable class of nonlinear partial differential equation systems }

For a positive viscosity constant $\nu >0$ consider the Cauchy problem 
\begin{equation}\label{cp1}
v_{i,t}-\nu \Delta v_i+ f_i(v,\nabla v)=0,~1\leq i\leq D,
\end{equation}
where $v=(v_1,\cdots,v_D)^T$ and $\nabla v=(\nabla v_1,\cdots ,\nabla v_D)^T$ are defined on the domain $[0,\infty)\times {\mathbb R}^D$ with data
\begin{equation}\label{cp2}
v_i(0,.)=h_i,~1\leq i\leq D.
\end{equation}
Here, $D\geq 3$ denotes the dimension of the problem and ${\mathbb R}$ denotes the field of real numbers.
\begin{rem}
The operators $f_i$ in (\ref{cp1}) can be global operators, e.g. partial integro-differential equation as in the case of the(Leray projection form of the) incompressible Navier Stokes equation. We shall define a scheme where  a solution is bounded in a ball of finite radius $C>0$ with respect to a regular norm, and such that this bound is preserved in time, i.e., there is a finite constant $C>0$  such that for $m\geq 2$, all $t_0\geq 0$ and some $\Delta >0$ (independent of $t_0\geq 0$) we have
\begin{equation}
\begin{array}{ll}
\max_{1\leq i\leq D}{\big |}h_i{\big |}_{H^m\cap C^m}\leq C,
~\mbox{ and }~\\
\\
\max_{1\leq i\leq D}{\big |}v_i(t_0,.){\big |}_{H^m\cap C^m}\leq C\Rightarrow  \max_{1\leq i\leq D}{\big |}v_i(t_0+\Delta,.){\big |}_{H^m\cap C^m}\leq C.
\end{array}
\end{equation}
In general our construction implies that $C>\max_{1\leq i\leq D}{\big |}h_i{\big |}_{H^m\cap C^m}$, and a spatial scaling parameter indicates this deviation of a strong semigroup property.
\end{rem}

Assume that the following conditions are satisfied.
\begin{itemize}
\item[a)] For $m\geq 2$ and $1\leq i\leq D$ the data satisfy
\begin{equation}
h_i\in H^m\cap C^m.
\end{equation}
Here, $C^m$ is the function space of continuous functions with continuous multivariate derivatives up to order $m$, and $H^m$ is the standard Sobolev space of order $m$, i.e., the space of functions with multivariate weak derivatives in $L^2$ up to order $m$. For $m=0$ the function space $C:=C^0$ is the set of continuous functions. We shall use the norm
\begin{equation}
{\big |}h_i{\big |}_{H^m\cap C^m}=\sum_{0\leq |\alpha|\leq m}{\big |}D^{\alpha}_xh_i{\big |}_{L^2\cap C},
\end{equation}
where for a continuous and bounded function $g:{\mathbb R}^D\rightarrow {\mathbb R}$ we define 
\begin{equation}
{\big |}g{\big |}_{L^2\cap C}:=\max\left\lbrace {\big |}g{\big |}_{L^2},\sup_{x\in {\mathbb R}^n}{\big |}g(x){\big |}\right\rbrace. 
\end{equation}
In some situations, it is possible to transform to problems on compact domains.   
If $\Omega \subset {\mathbb R}^n$ is compact, then ${\big |}g{\big |}_{L^2(\Omega)\cap C(\Omega)}$ denotes a local form of this norm. 

 \item[b)] The nonlinear terms $f_i$ satisfy a sub-homogeneity condition in the sense that for spatial scaling $y=rx,~r>0$ (where $r$ is any given positive real number), and for functions $v_i,~1\leq i\leq D$ with $v_i(t,.)\in H^m\cap C^m,~m\geq 2$, we have for all $1\leq i\leq D$, all $t\geq 0$, and all $y=rx\in {\mathbb R}^D$
\begin{equation}\label{fiin}
{\big |}f_i(v(t,x),\nabla v(t,x)){\big |}\leq r{\Big |}f_i(v(t,x),\frac{1}{r}\nabla v(t,x)){\Big |}.
\end{equation}
Note that for the transformation $v^r_i(t,y)=v_i(t,x)$ we have 
\begin{equation}
r f_i\left( v(t,x),\frac{1}{r}\nabla v(t,x)\right) =r f_i(v^r(t,y),\nabla v^r(t,y)),
\end{equation}
because $v_{i,j}=v^{r}_{i,j}\frac{dy_j}{dx_j}=v^{r}_{i,j}r$. Hence, the inequality in (\ref{fiin}) may be also rewritten such that for any given $r>0$ and with $v^r_i(t,y)=v_i(t,x)$ we have 
\begin{equation}\label{fiin}
{\big |}f_i(v(t,x),\nabla v(t,x)){\big |}\leq r{\Big |}f_i(v^r(t,y),\nabla v^r(t,y)){\Big |}.
\end{equation}
\item[c)] 
In the following for functions  $F^i_j$ we use the notation $$F^i_{j,j}(v,\nabla v):=\frac{\partial}{\partial x_j}F^i_{j}(v(t,x),\nabla v(t,x)),$$ i.e., the notation refers to derivatives with respect to the argument $x_j$ of the composition; some readers may prefer the notation $$\frac{\partial}{\partial x_j}F^i_{j}(v(t,x),\nabla v(t,x))=\left( F^i_{j}(v(t,x),\nabla v(t,x))\right)_{,j}, $$ which we avoid for simplicity of notation. Multivariate derivatives with respect to the arguments of $F^i_j$ may be denoted by $D^{\gamma}_{v,\nabla v} F^i_j$, where in this case $\gamma$ is a multiindex of length $D+D^2$. Note that these are usual multivariate derivatives.
We assume that there exists a matrix of regular functions $(F^i_j)_{1\leq i,j\leq D}$ such that for all $1\leq i\leq D$ 
\begin{equation}\label{Fijf}
\sum_{j=1}^D F^i_{j,j}(v,\nabla v)=f_i(v,\nabla v).
\end{equation}
\begin{rem}
The assumption of a matrix $F^i_j$ with (\ref{Fijf}) implies that we have classical local time representations in terms of the {\it first spatial derivatives} of the Gaussian. Symmetries of the first order spatial derivatives of the Gaussian can be exploited if the convoluted term in the local representation satisfies a Lipschitz condition. This structure together with the scaling assumption of item b) ensure that local diffusion damping can offset possible  growth of the nonlinear terms for a spatially scaled model. 
\end{rem}
We assume that $f_i$ satisfies a local Lipschitz condition, i.e., for $g=(g_1,\cdots,g_D)$,~$g_i\in C^m\cap H^m$ and for any finite constant $C>0$  there exists a finite constant $L$ (depending on $C$) such that
\begin{equation}
{\big |}g_i{\big |}_{H^m\cap C^m}\leq C\Rightarrow  {\big |}f_i(g,\nabla g){\big |}_{H^{m-1}\cap C^{m-1}}\leq L {\big |} g{\big |}_{H^{m-1}\cap C^{m-1}}.
\end{equation}
Note the loss of one order of regularity in this stipulation.
Furthermore, if $F^i_j$ has a local interpretation 
then we assume that for all multiindices $0\leq |\beta|\leq m$ the functions $D^{\beta}_{v,\nabla v}F^i_{j}:{\mathbb R}^{D+D^2}\rightarrow {\mathbb R}$ are locally Lipschitz continuous (on any compact domain in $D\subset {\mathbb R}^D$ and with respect to all arguments).
%
In any case we  require that for any $g=(g_1,\cdots,g_D)$ with ${\big |}g_i{\big |}_{H^m\cap C^m}\leq C <\infty $ we have Lipschitz continuity of $ y\rightarrow D^{\gamma}_xF^i_j\left(g(y),\nabla g(y) \right)$ with finite Lipschitz constants $L^{i\gamma}_j $ (dependent on $C$) in the sense that for all $1\leq i,j\leq D$ and all $0\leq |\gamma|\leq m-1$ and $y,y'\in {\mathbb R}^D$
\begin{equation}
{\big |}F^i_j\left(g(y),\nabla g(y) \right)-F^i_j\left(g(y'),\nabla g(y') \right){\big |}\leq L^{i0}_j|y-y'|,
\end{equation}
and
\begin{equation}
{\big |}D^{\gamma}_xF^i_{j,j}\left(g(y),\nabla g(y) \right)-D^{\gamma}_xF^i_{j,j}\left(g(y'),\nabla g(y') \right){\big |}\leq L^{i\gamma}_j|y-y'|.
\end{equation}
Here, $D^{\gamma}_x$ denotes the multivariate spatial derivative with respect to the multiindex $\gamma=(\gamma_1,\cdots,\gamma_D)$ and with respect to the argument of $g$.
\item[d)] The verification of the technical condition in c) can be simplified for specific models if we add a stronger assumption concerning the data. This additional assumption is also useful, if generalisations of the diffusion term or viscosity limits are considered, e.g., if we replace the Laplacian term by a H\"{o}rmander vector-field condition for highly degenerate operators of second order. Let
\begin{equation}
h_i\in H^m\cap C^m\cap  {\cal C}^{m(D+1)}_{pol,m},
\end{equation}
where for a given integer $l\geq 1$
\begin{equation}\label{pol1}
{\cal C}^{l}_{pol,m}={\Big \{} f:{\mathbb R}^D\rightarrow {\mathbb R}: \\
\\
\exists c>0~\forall |x|\geq 1~\forall 0\leq |\gamma|\leq m~{\big |}D^{\gamma}_xf(x){\big |}\leq \frac{c}{1+|x|^l} {\Big \}}.
\end{equation}
Note that the latter function space has a multiplicative property, i.e., $g,h\in {\cal C}^{l}_{pol,m}$ implies that $gh\in {\cal C}^{2l}_{pol,m}$.
In addition to (\ref{pol1}) we then require a 'submultiplicative property of order $k\in \left\lbrace 0,\cdots,m(D+1)-1\right\rbrace$' , i.e., for $g=(g_1,\cdots,g_D)$ along with $g_i\in {\cal C}^{m(D+1)}_{pol,m}$ we require that for some $k\in \left\lbrace 0,\cdots,m(D+1)-1\right\rbrace$ we have
\begin{equation}
g_i\in {\cal C}^{m(D+1)}_{pol,m}, (\nabla g_i)_j\in {\cal C}^{m(D+1)}_{pol,m} \mbox{ implies } D^{\beta}_xF^i_j(g,\nabla g)\in {\cal C}^{2m(D+1)-k}_{pol,m}. 
\end{equation}
Even the latter condition is a sufficient for a class of operators which includes some models of fluid mechanics. Especially, this 'multiplicative property' of the nonlinear term holds for the incompressible Navier-Stokes equation operator. The number $m(D+1)$ in the upperscript of ${\cal C}^{m(D+1)}_{pol,m}$ is not a sharp choice, but a choice which may be also sufficient for the consideration of viscosity limits. 
\end{itemize}

We have
\begin{thm}\label{mainthm}
If the set of conditions a), b), c)  or the stronger set of conditions a), b), c) and d) are satisfied, then the Cauchy problem in (\ref{cp1}), (\ref{cp2}) has a global classical solution $v_i\in C^1\left(\left[0,\infty\right) , H^m\cap C^m \right),~1\leq i\leq D$. 
\end{thm}

Some remarks are in order.
\begin{rem}
We have described a quadratic system where a vector $v$ with $D$ components solves $D$ equations. This is not an essential restriction, and it is convenient.
\end{rem}
\begin{rem}
We do not claim uniqueness in (\ref{mainthm}) although in many special situations standard arguments may lead to uniqueness. Note that in some situations the method described below may be applied in order to get global solution branches of viscosity limits of the equations in (\ref{cp1}), and then there are examples, where determinism is lost and global regular solution branches exist next to singular solutions.
\end{rem}
\begin{rem}
We use Gaussian upper bound estimate. We cannot prove a strong contractive semi-group property for a class of nonlinear operators which includes the Navier Stokes equation operator. The deviation from a strong semigroup property is indicated directly, if we estimate for $v^r_i,~1\leq i\leq D$ with $r>1$. In any case estimates are obtained for multiples of a discrete time, and then the upper bound for all time is constructed by a local time contraction result in regular space, which gives an upper bound in regular space, and this upper bound is larger than the data in the same regular space in general.
\end{rem}
\begin{rem}
We may reduce to a local Lipschitz continuity assumption in c) if we add d) or the weaker assumption that for $u=(u_1,\cdots,u_D)$ with $u_i\in H^m\cap C^m$ for $m\geq 2$ we have that $F^i_j(u,\nabla u)\in L^2$. 
\end{rem}

\begin{exmp}\label{exmpns}
In case of the incompressible Navier Stokes equation (cf. \cite{LL} for the modeling) the incompressibility condition implies
\begin{equation}
\sum_{j=1}^D\frac{\partial ( v_iv_j)}{\partial x_j}=\sum_{j=1}^Dv_j\frac{\partial v_i}{\partial x_j}+v_i\sum_{j=1}^D\frac{\partial v_j}{\partial x_j}= \sum_{j=1}^Dv_j\frac{\partial v_i}{\partial x_j}.
\end{equation}
Therefore we may define
\begin{equation}\label{Fij}
F^i_j=v_jv_i-\delta_{ij}K_D \ast_{sp}\sum_{l,m=1}^Dv_{l,m}v_{m,l},
\end{equation}
where $K_D$ denotes the Laplacian kernel of dimension $D$ and $\delta_{ij}$ is the Kronecker-$\delta$. Here $\ast_{sp}$ denotes the convolution with respect to the spatial variables. 
Note that for $y=rx$, $v^r_i(t,y)=v_i(t,x)$, and $z=rw$, and the first spatial derivative of the Laplacian kernel $x\rightarrow K_{D,i}(x)=\frac{x_i}{|x|^D}$ in the Leray projection term (written as an operator on the Jacobian $J(v)=\left(v_{l,m} \right)_{1\leq l,m\leq D}$) satisfies
\begin{equation}
\begin{array}{ll}
L(J(v))=\int_{{\mathbb R}^D}\frac{x_i-w_i}{|x-w|^D}\sum_{l,m=1}^Dv_{l,m}(t,w)v_{m,l}(t,w)dw\\
\\
=\int_{{\mathbb R}^D}
\frac{\frac{y_i-z_i}{r}}{\frac{|y-z|^D}{r^D}}
\sum_{l,m=1}^Drv^r_{l,m}(t,z)rv^r_{m,l}(t,z)\frac{1}{r^D}dz\\
\\
\int_{{\mathbb R}^D}
\frac{(y_i-z_i)r^{D-1}}{|y-z|^D}
\sum_{l,m=1}^Drv^r_{l,m}(t,z)rv^r_{m,l}(t,z)\frac{1}{r^D}dz\sim rL(J(v^r)).
\end{array}
\end{equation}
This linear spatial scaling of the Leray projection term can be observed also from the linear scaling of the pressure gradient. The Poisson equation for the scaled equation becomes $r^2\sum_{l,m=1}^D v^r_{l,m}v^r_{m,l}=r^2 \Delta p^r$ and is identical to the pressure elimination equation in original coordinates as the factor $r^2$ cancels. More explicitly, for $y=rx$ and $p^r(y)=p(x)$ we have the spatial derivative transformation $(\nabla)_ip:=p_{,i}=\frac{\partial p}{\partial x_i}=\frac{\partial p^r}{\partial y_i}\frac{dy_i}{dx_i}=p^r_{,i}r$, and the original Navier Stokes equation transforms to
\begin{equation}\label{nsr1}
\frac{\partial v^r_i}{\partial t}-r^2\nu\Delta v^r_i+r\sum_{j=1}^Dv^r_j\frac{\partial v^r}{\partial x_j}=r(\nabla)_i p^r,
\end{equation}
and, applying the divergence operator and incompressibility, i.e., $\div v=0$, leads indeed to
\begin{equation}
r^2\sum_{l,m=1}^D v^r_{l,m}v^r_{m,l}=r^2 \Delta p^r\Leftrightarrow \sum_{l,m=1}^D v^r_{l,m}v^r_{m,l}= \Delta p^r.
\end{equation}
Hence, the Leray projection form of the scaled equation becomes
\begin{equation}\label{nsr2}
\frac{\partial v^r_i}{\partial \tau}-r^2\nu\Delta v^r_i+r\sum_{j=1}^Dv^r_j\frac{\partial v^r}{\partial x_j}=r\int_{{\mathbb R}^D}K_{D,i}(.-y)\sum_{l,m=1}^D v^r_{l,m}v^r_{m,l}(y)dy
\end{equation}
If the additional assumption d) holds, then it is rather obvious that the technical Lipschitz condition in c) holds, but it can also be verified under the weaker set of assumptions a), b) and c). 
Assume that at time $t_0\geq 0$ we have the upper bound
\begin{equation}
\max_{1\leq i\leq D}{\big |}v^r_i(t_0,.){\big |}_{H^m\cap C^m}\leq C \mbox{ for some $m\geq 2$,}
\end{equation}
for some finite constant $C>0$. Local time iteration schemes with respect strong norms lead to local time representations
\begin{equation}\label{rep1}
\begin{array}{ll}
v^r_i=v^r_i(t_0,.)\ast_{sp} G^r_{\nu}+r\sum_{j=1}^Dv^r_j\frac{\partial v^r}{\partial x_j}\ast G^r_{\nu}\\
\\
+r\int_{{\mathbb R}^D}K_{D,i}(.-y)\sum_{l,m=1}^D v^r_{l,m}(.,y)v^r_{m,l}(.,y)dy\ast G^r_{\nu},
\end{array}
\end{equation}
on a short time interval $\left[t_0,t_0+\Delta \right]$ for small $\Delta >0$, and where $\ast$ denotes convolution with respect to space and time on this time interval. Here, $G^r_{\nu}$  is the fundamental solution of $q_{,t}-r^2\nu \Delta q=0$ (considered on the time interval $\left[t_0,t_0+\Delta \right]$). Similar local time representations hold for spatial derivatives of velocity function. Indeed for multiindices $\beta=(\beta_1,\cdots,\beta_D)$ and $\gamma=(\gamma_1,\cdots,\gamma_d)$  and $1\leq |\beta|=1+|\gamma|\leq m$ with $\beta_p=\gamma_p+1$ and $\beta_j=\gamma_j$ for $p\neq j$ we have on the time interval $[t_0,t_0+\Delta]$
\begin{equation}\label{rep2der}
\begin{array}{ll}
D^{\beta}_xv^r_i=D^{\beta}_xv^r_i(t_0,.)\ast_{sp} G^r_{\nu}+rD^{\gamma}_x\left( \sum_{j=1}^Dv^r_j\frac{\partial v^r}{\partial x_j}\right) \ast G^r_{\nu,p}\\
\\
+rD^{\gamma}_x\left( \int_{{\mathbb R}^D}K_{D,i}(y)\sum_{l,m=1}^D v^r_{l,m}(.,x-y)v^r_{m,l}(.,x-y)dy\right) \ast G^r_{\nu,p}.
\end{array}
\end{equation}
Local time contraction shows that 
\begin{equation}
\max_{1\leq i\leq D}\sup_{t\in [t_0,t_0+\Delta]}{\big |}v^r_i(t,.){\big |}_{H^m\cap C^m}\leq C+1
\end{equation}
such that we have the Leray projection estimate 
\begin{equation}\label{repleray}
\begin{array}{ll}
r{\Big |} \left( \left( K_{D,i}\ast_{sp}\sum_{l,k=1}^D D^{\gamma}_x\left(v^r_{l,k}v^r_{k,l}\right) \right) \ast G^r_{\nu,p}\right) (t,.){\Big |}_{L^2\cap C}\\
\\
\leq rc^D_m (C+1) \max_{l=1}^D\max_{|\delta|\leq m-1} {\Big |} \left( \left( K_{D,i}\ast_{sp}
D^{\delta}_xv^r_{l} \right) \ast 
G^r_{\nu,p}\right) (t,.){\Big |}_{L^2\cap C},
\end{array}
\end{equation}
where $c^D_m$ is a finite constant depending only on dimension $D$ and the regularity number $m$. Since for $t\in [t_0,t_0+\Delta]$ and $0\leq |\delta|\leq m-1$
\begin{equation}
{\big |}D^{\delta}_xv^r_{l}(t,.){\big |}_{H^1\cap C^1}\leq C+1
\end{equation}
we have spatial Lipschitz continuity of $K_{D,i}\ast_{sp}D^{\delta}_xv^r_{l}(t,.)$ for multivariate spatial derivatives of order $|\delta|\leq m-1$ and in case of dimension $D=3$. As in the main theorem we then may use  symmetry of the first order derivative of the Gaussian (cf. below).  
Indeed, application of Theorem \ref{mainthm} implies the existence of a global regular  solution branch. In this special case uniqueness is implied. 
\end{exmp}

\section{Proof Theorem \ref{mainthm}}
\begin{itemize}
\item[i)] We rewrite the equation using the integrability of the nonlinear terms. From (\ref{cp1}) and assumption c) we have 
\begin{equation}
v_{i,t}-\nu \Delta v_i+ \sum_{j=1}^DF^i_{j,j}(v,\nabla v)=0,~1\leq i\leq D.
\end{equation}
Recall the notation $F^i_{j,j}(v,\nabla v)=\left(F^i_{j}(v,\nabla v) \right)_{,j}$, where the derivative is with respect to the $j$th spatial argument $x_j$ of the composition of $F^i_j$ with $(v,\nabla v)$.

\item[ii)] Assume that data $v_i(t_0,.),~1\leq i\leq D$ are given. We have local time contraction in spatial $H^m\cap C^m$ space. More precisely, define a time local iteration scheme $v^k_i,~1\leq i\leq D,~k\geq 0$ in a time interval $[t_0,t_0+\Delta]$, where $t_0\geq 0$ and where for $t\in [t_0,t_0+\Delta]$ we define
\begin{equation}
v^0_i(t_0,.)=v_i(t_0,.)\ast_{sp} G^{r}_{\nu},~\mbox{(n.b. $v^0_i(0,.)=h_i(.)$)}.
\end{equation}
Furthermore, for $k\geq 1$ $v^k_i$ is a solution of the linearized equation
\begin{equation}
v^k_{i,t}-\nu \Delta v^k_i+ \sum_{j=1}^DF^i_{j,j}(v^{k-1},\nabla v^{k-1})=0,~1\leq i\leq D,
\end{equation}
where $v^k_i(t_0,.)=v_i(t_0,.)$ for $1\leq i\leq D$.
Note that 
\begin{equation}
v^k_{i}=v_i(t_0,.)\ast_{sp}G^{r}_{\nu}+\sum_{j=1}^DF^i_{j,j}(v^{k-1},\nabla v^{k-1})\ast G^{r}_{\nu}
\end{equation}
for $1\leq i\leq D$.
\begin{lem}\label{contrlem}
Let $t_0\geq 0$ and assume that for some $m\geq 2$ and a finite constant $C_0>0$ we have
\begin{equation}
{\big |}v_i(t_0,.){\big |}_{H^m\cap C^m}\leq C_0.
\end{equation}
Then  there exists a $\Delta >0$ dependent only on dimension and on the constant ${\big |}v_i(t_0,.){\big |}_{H^m\cap C^m}$ such that on the time interval $[t_0,t_0+\Delta]$ the  functional increments  $\delta v^{k+1}_j=v^{k+1}_j-v^{k}_j,~1\leq j\leq D$ satisfy for $k\geq 1$
\begin{equation}
\sup_{t\in [t_0,t_0+\Delta]}{\big |}\delta v^{k+1}_j(t,.){\big |}_{H^m\cap C^m}\leq \frac{1}{2} \sup_{t\in [t_0,t_0+\Delta]}{\big |}\delta v^{k}_j(t,.){\big |}_{H^m\cap C^m}
\end{equation}
and
\begin{equation}
\sup_{\tau\in [t_0,t_0+\Delta]}{\big |}\delta v^{1}_i(t,.){\big |}_{H^m\cap C^m}\leq \frac{1}{2}.
\end{equation}

\end{lem}
The proof uses classical representations and the Lipschitz continuity of (spatial derivatives of) $F^i_j$ and is given in the appendix.

\item[iii)] The local time contraction result of item ii) implies that there exist regular local time solutions $v_i\in C^1\left(\left[t_0,t_0+\Delta \right],H^m\cap C^m  \right),~1\leq i\leq D$. 
\begin{rem}
If the $f_i$ satisfy the condition in d), then for the functional increments $\delta v^k=v^k-v^{k-1},~k\geq 1,$ we have $$\delta v^k_i\in C^1\left(\left[t_0,t_0+\Delta \right],H^m\cap C^m \cap {\cal C}^{m(D+1)}_{pol,m}  \right),~1\leq i\leq D,$$ and this holds also in the limit for $v_i,~1\leq i\leq D$. Here, note that 
\begin{equation}
\delta v^k_{i}=\sum_{j=1}^DF^i_{j,j}(v^{k-1},\nabla v^{k-1})\ast G^{r}_{\nu}-\sum_{j=1}^DF^i_{j,j}(v^{k-2},\nabla v^{k-2})\ast G^{r}_{\nu}
\end{equation}
for $1\leq i\leq D$, where the submuliplicative property of $f_i$ implies a preservation of the order of the order of spatial polynomial decay (as $k$ increases) where it clearly offsets possible effects of decrease the order of spatial polynomial decay caused
by convolutions with a first order spatial derivative of a Gaussian. For more complicated second order diffusions we have a convolution with a density which can cause a stronger decrease of polynomial order decrease at spatial infinity, and we need a stronger submultiplicative property (smaller $k$ in assumption  d)).
\end{rem}      
We consider local classical representations of solutions and use the convolution rule in order to get for all $t\in [t_0,t_0+\Delta]$ and $x\in {\mathbb R}^D$
\begin{equation}\label{locrep}
\begin{array}{ll}
v_i(t,x)=v_i(t_0,.)\ast_{sp}G_{\nu}+\sum_{j=1}^D F^i_{j,j}\left(v,\nabla v \right)\ast G_{\nu}\\
\\
=v_i(t_0,.)\ast_{sp}G_{\nu}+\sum_{j=1}^D F^i_{j}\left(v,\nabla v \right)\ast G_{\nu,j} 
\end{array}
\end{equation}
Again remember that $,j$ refers to derivative with respect to $x_j$, where we use the notation $F^i_{j,j}\left(v,\nabla v \right)=\left( F^i_{j,j}\left(v,\nabla v \right)\right)_{,j}$.
Furthermore, the symbol $\ast$ denotes convolution with respect to space and time.
In the last line the nonlinear terms are convoluted with first order spatial derivatives of the Gaussian, while the first term on the right side of (\ref{locrep}) is a convolution with the Gaussian, which behaves completely different. Here we may use only one spatial scaling parameter $r>0$ and exploit spatial effects of the operator. We shall observe that small damping of the latter term offsets possible growth caused by the former.
Note that for spatial multivariate derivatives of order $|\beta|$ and for $0\leq |\gamma|+1=|\beta|\leq m$, $\beta_k=\gamma_k+1$, and $\beta_l=\gamma_l$ for $l\neq k$ we have time local representations of the form
\begin{equation}\label{locrepaa}
\begin{array}{ll}
D^{\beta}_xv_i(t,x)=D^{\beta}_xv_i(t_0,.)\ast_{sp}G_{\nu}+\sum_{j=1}^D D^{\gamma}_xF^i_{j,j}\left(v,\nabla v \right)\ast G_{\nu,k}.
\end{array}
\end{equation}

\item[iv)] Using Lipschitz estimates we get an upper bound of the nonlinear terms.
First note that first order spatial derivatives of the Gaussian $G_{\nu,j}$ have a symmetry which can be combined with Lipschitz continuity of the nonlinear function terms. We have
\begin{equation}
G_{\nu,j}(t,y)=\frac{-y_j}{2\nu t}\frac{1}{\sqrt{4\pi \nu t}^D}\exp\left(-\frac{|y|^2}{4 \nu t} \right). 
\end{equation}
Define $y^{j,-}=(y^{j,-}_1,\cdots,y^{j,-}_D)$, where $y^{j,-}_k=y_k$ for $k\neq j$ and $y^{j,-}_j=-y_j$.
Assuming ${\big |}v_i(t_0,.){\big |}_{H^m\cap C^m}=C_{t_0}<\infty$ the local contraction results tells us that for some $\Delta >0$ and all time $t\in [t_0,t_0+\Delta]$  we have 
\begin{equation}
{\big |}v_i(t,.){\big |}_{C^m\cap H^m}\leq C_{t_0}+1.
\end{equation}
Recall that functions $F^i_j:{\mathbb R}^{D+D^2}\rightarrow {\mathbb R}$ are Lipschitz on any finite ball, especially on a ball of $B_{C_{t_0}+1}(0)$   of radius $C_{t_0}+1$ around zero.
It follows that for all $1\leq i,j\leq D$
\begin{equation}\label{lip1}
\begin{array}{ll}
{\Big |}\int_{{\mathbb R}^D}F^i_j(v,\nabla v)(t,x-y)\frac{-y_i}{2\nu t}\frac{1}{\sqrt{4\pi \nu t}^D}\exp\left(-\frac{|y|^2}{4 \nu t} \right)dy {\Big |}\\
\\
={\Big |}\left( \int_{{\mathbb R}^D}\left( F^i_j(v,\nabla v)(t,x-y)-F^i_j(v,\nabla v)(t,x-y^{j,-})\right)\right)\times\\
\\
 \times\frac{|y_j|}{2 \nu t}\frac{1}{\sqrt{4\pi \nu t}^D}\exp\left(-\frac{|y|^2}{4 \nu t} \right)dy {\Big |}\\
\\
\leq 2L^{i0}_j{\Big |}\int_{{\mathbb R}^D}\frac{y^2_j}{\sqrt{4\pi \nu t}^D}\exp\left(-\frac{|y|^2}{4 \nu t} \right)dy {\Big |},
\end{array}
\end{equation}
where we use assumption c). Analogous estimates hold for spatial derivatives with the related Lipschitz constants $L^{i\gamma}_j$.

\item[v)] The nonlinear upper bound of item iv) has a scaling which is different from a normal Gaussian. Moreover the linear and the nonlinear part of the equation have scaling constraints due to the sub-homogeneity condition.
For $ t=\tau\rho$ and $rx=y$, and $v^{\rho,r}_i(\tau,y)=v_i(t,x)$ we have $v_{i,t}=v^{\rho,r}_{i,\tau}\frac{d\tau}{d t}=v^{\rho,r}_{i,\tau}\frac{1}{\rho}$,  $v_{i,j}=v^{\rho,r}_{i,j}r$, and $v_{i,j,j}=v^{\rho,r}_{i,j,j}r^2$. 
\begin{equation}\label{cp1r}
v^{\rho,r}_{i,\tau}-\rho r^2\nu \Delta v^{\rho,r}_i+ \rho f_i(v^{\rho,r},r\nabla v^{\rho,r})=0,~1\leq i\leq D.
\end{equation}
Now,
\begin{equation}
\begin{array}{ll}
{\big |}f_i(v^{\rho,r}(t,y),r\nabla v^{\rho,r}(t,y)){\big |}={\big |}f_i(v(t,x),\nabla v(t,x)){\big |}\\
\\
\leq r{\big |}f_i(v(t,x),\frac{1}{r}\nabla v(t,x)){\big |}=r{\big |}f_i(v^r(t,y),\nabla v^r(t,y)){\big |}.
\end{array}
\end{equation}

Hence,
\begin{equation}
v^{\rho,r}_{i,\tau}-\rho r^2\nu \Delta v^{\rho,r}_i+ \rho\sum_{j=1}^DF^i_{j,j}(v^{\rho,r},r\nabla v^{\rho,r})=0,~1\leq i\leq D,
\end{equation}
where
\begin{equation}
\begin{array}{ll}
{\big |}\rho\sum_{j=1}^DF^i_{j,j}(v^{\rho,r},r\nabla v^{\rho,r}){\big |}={\big |}\rho f_i(v^{\rho,r},r\nabla v^{\rho,r}){\big |}\\
\\
\leq \rho r{\big |} f_i(v^{\rho,r},\nabla v^{\rho,r}){\big |}=\rho r{\big |}\sum_{j=1}^DF^i_{j,j}(v^{\rho,r},\nabla v^{\rho,r}){\big |}.
\end{array}
\end{equation}
The local contraction result transfers to the scaled situation. It follows that on a time interval $[t_0,t_0+\Delta_0]$ with $\Delta_0=\rho\Delta$ we have
\begin{equation}
\sup_{t\in [t_0,t_0+\Delta_0]}{\big |}v_i(t,.){\big |}_{H^m\cap C^m}{\big |}v_i(t_0,.){\big |}_{H^m\cap C^m}+1:=C_{t_0m}+1.
\end{equation}
Note that this holds for $t_0=0$ especially.

For the scaled function we have the local representation
\begin{equation}\label{locreplat}
\begin{array}{ll}
v^{\rho,r}_i(\tau,x)=v^{\rho,r}_i(t_0,.)\ast_{sp}G^{\rho,r}_{\nu}+\rho r\sum_{j=1}^D F^i_{j,j}\left(v^{\rho,r},\nabla v^{\rho,r} \right)\ast G^{\rho,r}_{\nu}\\
\\
=v^{\rho,r}_i(t_0,.)\ast_{sp}G^{\rho,r}_{\nu}+\rho r\sum_{j=1}^D F^{i}_{j}\left(v^{\rho,r},\nabla v^{q,r} \right)\ast G^{\rho,r}_{\nu,j} 
\end{array}
\end{equation}
For the scaled Gaussian $G^{\rho,r}_{\nu,i}(\tau,x):=G_{\nu}(t,y)$ we have
\begin{equation}
\begin{array}{ll}
{\big |}G^{\rho,r}_{\nu,i}(\tau,y){\big |} 
={\Big |}\frac{-2y_i}{4\pi \rho r^2\nu \tau}\frac{1}{\sqrt{4\pi\rho r^2\nu t}^D}\exp\left(-\frac{|y|^2}{4\rho r^2\nu \tau} \right){\Big |}\\
\\
\leq \frac{2}{(4\pi\nu\rho r^2 \tau)^{\delta}|y|}\left( |y|^2\right)^{\delta-D/2} \left( \frac{|y|^2}{4\pi\rho r^2\nu \tau}\right)^{D/2+1-\delta} \exp\left(-\frac{|y|^2}{4\rho r^2\nu \tau} \right).
\end{array}
\end{equation}
Hence we have for $\delta \in (0,1)$ and all $\rho,r>0$ 
\begin{equation}\label{grrh0}
{\big |}G^{\rho,r}_{\nu ,i}(\tau,y){\big |}\leq \frac{C}{ (4\pi \rho r^2 \nu\tau)^{\delta}|y|^{D+1-2\delta}}, 
\end{equation}
where the upper bound constant $$C=\sup_{|z|>0}\left( z\right)^{D/2+1-\delta} \exp\left(-z^2\right) >0$$ is sufficient and independent of $\nu >0$. 
We are interesting in the scaling of convolutions of these Gaussians with (globally) Lipschitz continuous functions $y\rightarrow l(y)$ with upper bound $l_0|y|$. Now, as $\nu >0$ we may choose the spatial parameter $r$ to be large enough such that
\begin{equation}
4\pi \rho r^2 \nu\geq 1.
\end{equation}

Using the abbreviations $\sigma_{\tau}=\tau -\sigma$, and $B=\left\lbrace y{\big |}|y|\leq  4\pi \rho r^2 \nu\right\rbrace$, we find
\begin{equation}\label{Best}
\begin{array}{ll}
\int_{t_0}^{\tau}\int_{B}l_0|y|{\big |}G^{\rho,r}_{\nu ,i}(\sigma_{\tau},y){\big |}dyd\sigma \leq \\
\\
\int_{t_0}^{\tau}\int_{B}\frac{l_0}{(4\pi\nu\rho r^2 \sigma_{\tau})^{\delta}}\left( |y|^2\right)^{\delta-D/2} \left( \frac{|y|^2}{4\pi\rho r^2\nu \sigma_{\tau}}\right)^{D/2+1-\delta} \exp\left(-\frac{|y|^2}{4\rho r^2\nu \sigma_{\tau}} \right)dyd\sigma\\
\\
\leq \int_{t_0}^{\tau}\int_{B}\frac{l_0}{(4\pi\nu\rho r^2 \sigma_{\tau})^{\delta}}\left( |y|^2\right)^{\delta-D/2} Cdyd\sigma\leq l_0C^*(\tau-t_0)^{1-\delta}
(4\pi\nu\rho r^2 )^{\delta},
\end{array}
\end{equation}
for any $\delta\in (0,1)$ and for some finite constant $C^*$ which is independent of the parameters $\rho,r,\nu$ (radial coordinates with measure $r^{D-1}dr$ may be used for the latter observation). Note that for small time and/or small parameters the corresponding integral on the complementary domain $\int_{t_0}^{\tau}\int_{{\mathbb R}^D\setminus B}\cdots$ becomes small such that (\ref{Best}) is indeed the essential estimate in the sense that it gives the whole upper bound up to an 'small $\epsilon$'. 
Now we can estimate the nonlinear term rephrasing (\ref{lip1}) above and including time. Using the sub- homogeneity property we have the upper bound 
\begin{equation}\label{lip2}
\begin{array}{ll}
{\Big |}\sum_{j=1}^D\rho r\int_{t_0}^{\tau}\int_{{\mathbb R}^D}F^i_j(v^r,\nabla v^r)(s,y-z)\times\\
\\
\times \frac{-2z_i}{4\pi \nu (\tau-s)}\frac{1}{\sqrt{4\pi \nu (\tau-s)}^D}\exp\left(-\frac{|z|^2}{4\nu\rho r^2 (\tau-s)} \right)dz ds{\Big |}\\
\\
={\Big |}\sum_{j=1}^D\rho r\left( \int_{t_0}^{\tau}\int_{{\mathbb R}^D}\left( F^i_j(v^r,\nabla v^r)(s,y-z)-F^i_j(v^r,\nabla v^r)(s,y-z^{j,-})\right)\right)\times\\
\\
 \times\frac{|y_j|}{2\pi \nu s}\frac{1}{\sqrt{4\pi \nu (\tau- s)}^D}\exp\left(-\frac{|y|^2}{4\rho r^2\nu (\tau-s)} \right)dyds {\Big |}\\
\\
\leq \rho r\sum_{j=1}^D2L^{i0}_j{\Big |}\int_{t_0}^{\tau}\int_{{\mathbb R}^D}\frac{y^2_j}{\sqrt{4\pi \rho r^2\nu (\tau-s)}^D}\exp\left(-\frac{|y|^2}{4\rho r^2\nu (\tau-s)} \right)dy {\Big |}\\
\\
\leq \rho r\sum_{j=1}^D2L^{i0}_j(4\pi \rho r^2 \nu )^{\delta}(C^*(\tau-t_0)^{1-\delta}+\epsilon),
\end{array}
\end{equation}
where we have a finite constant $C^*$ with 
\begin{equation}\label{c*}
\begin{array}{ll}
C^*(\tau-t_0)^{1-\delta}+\epsilon=2D\int_{t_0}^{t_0+\Delta}\int_{B_1(0)}\frac{C}{ \sigma^{\delta}|y|^{D+1-2\delta}}dyd\sigma\\
\\
+\sum_{j=1}^D\int_{t_0}^{\tau}\int_{{\mathbb R}^D\setminus B_1(0)}|y_j||G^{\rho,r}_{\nu,j}(\sigma,y)|d\sigma dy.
\end{array}
\end{equation}
Here $\epsilon$ is defined by the second term on the right side of (\ref{c*}).
The first summand on the right side of (\ref{c*}) is independent of $r,\rho,\nu$.
The last summand on the right side of (\ref{c*}) is relatively small (goes to zero with exponential decay as $\rho r^2\nu$ or $\Delta_0$ go to zero). Hence $\epsilon$ is small compared to any of the polynomials of lower order which determine the main part of $C^*$. Summing up we have
\begin{equation}\label{lip2a}
\begin{array}{ll}
{\Big |}\sum_{j=1}^D\rho r\int_{t_0}^{\tau}\int_{{\mathbb R}^D}F^i_j(v^r,\nabla v^r)(s,y-z)\frac{-2z_i}{4 \nu s}\frac{1}{\sqrt{4\pi \nu (\tau-s)}^D}\exp\left(-\frac{|z|^2}{4 \nu\rho r^2 (\tau-s)} \right)dz {\Big |}\\
\\
\leq \rho rL_0(4\pi \rho r^2 \nu )^{\delta}((\tau-t_0)^{1-\delta}+\epsilon),
\end{array}
\end{equation}
where $L_0:=(\sum_{j=1}^D2L^{i0}_j)$ and with an $\epsilon >0$ of exponential decay with respect to $\nu\rho r^2\Delta_0$ as $\nu\rho r^2\Delta_0$ becomes small. Analogous considerations lead to an upper bound
\begin{equation}\label{lip2a}
\begin{array}{ll}
{\Big |}\sum_{j=1}^D\rho r\int_{t_0}^{\tau}\int_{{\mathbb R}^D}D^{\beta}_xF^i_j(v^r,\nabla v^r)(\tau-s,y-z)\frac{-z_i}{4\pi \nu s}\frac{1}{\sqrt{4\pi \nu t}^D}\exp\left(-\frac{|z|^2}{4 \nu\rho r^2 t} \right)dz {\Big |}\\
\\
\leq \rho rL_m(4\pi \rho r^2 \nu )^{\delta}((\tau-t_0)^{1-\delta}C^*+\epsilon),
\end{array}
\end{equation}
for some constant $L_m$ which is independent of $\rho,r$ and $\nu$ and which we choose such that it serves for all spatial derivatives up to order $m$. We may assume that $L_m\geq 0$ such that we may use $L_m$ for all these upper bounds as a constant which is independent of specific multiindices less or equal to $m$.

\item[vi)] We use small damping of the heat convolution on infinite space, if a threshold is exceeded, i.e., we may assume that for the data we have for all $0\leq |\beta|\leq m$
\begin{equation}\label{ass}
\max_{1\leq i\leq D}{\big |}D^{\beta}_xv_i(t_0,.){\big |}_{L^2\cap C}\geq 1.
\end{equation}
If the latter condition is not satisfied for some $\beta$ then there is a $\Delta>0$ such that the respected norm is less or equal to $1$ for some time $t\in [t_0,t_0+\Delta]$, and we need no damping estimate for this part of the $H^m\cap C^m$-norm in the interval $[t_0,t_0+\Delta]$. This way we construct an upper bound close to a constant $C_m$, where $C_m$ is $\max_{1\leq i\leq D}{\big |}h_i{\big |}_{H^m\cap C^m}$ plus the number of terms in the standard definition of the $H^m\cap C^m$-norm. 
 We consider $L^2$-estimates. We apply a Fourier transform with respect to the spatial variables, i.e., the operation
\begin{equation}
{\cal F}(u)(\tau,\xi)=\int_{{\mathbb R}^D}\exp\left(-2\pi i x\xi\right)u(\tau,x)dx,
\end{equation}
in order to analyze the viscosity damping encoded in the first term on the right side of (\ref{locrep}) on a time interval $\left[t_0,t_0+\Delta_0 \right]$, where $\Delta_0=\rho\Delta$. For $\tau\in \left[t_0,t_0+\Delta_0 \right]$ and parameters $r,\rho >0$ we have
\begin{equation}
\begin{array}{ll}
{\cal F}\left( v^{\rho ,r,\nu}_i(t_0,.)\ast_{sp}G^{\rho,r}_{\nu}(\tau-t_0)\right)={\cal F}\left( v^{\rho,r,\nu}_i(t_0,.)\right) {\cal F}\left( G^{\rho,r}_{\nu}(\tau-t_0,.)\right)\\
\\
={\cal F}\left( v^{\rho,r,\nu}_i(t_0,.)\right) \exp\left(-4\pi^2\rho r^2\nu (\tau-t_0) (.)^2 \right),
\end{array}
\end{equation}
where we use (let $t_0=0$ for simplicity)
\begin{equation}
\begin{array}{ll}
{\cal F}\left( G^{\rho,r}_{\nu}(\tau,.)\right) (\tau,\xi)=
{\cal F}\left( \frac{1}{{\sqrt{4\pi \rho r^2 \nu \tau}^D}}
\exp\left(-\frac{(.)^2}{4\nu \rho r^2 \tau} \right)\right)  (\tau,\xi)\\
\\
=\exp\left(-4\pi^2\rho r^2\nu \tau |\xi|^2 \right).
\end{array}
\end{equation}
In the following we let $t_0=0$ and remark that the following estimates hold for $t_0 >0$ if $\tau$ is replaced by $\tau -t_0$.
For $\Delta >0$ small enough (such that, say, $8\pi^2\rho r^2\nu \tau \Delta_0^2\leq 1$), and for $\tau\in [0,\Delta_0]$ we get
\begin{equation}
\begin{array}{ll}
{\big |}v^{\rho,r,\nu}_i(t_0,.)\ast_{sp}G^{\rho,r}_{\nu}(\tau,.){\big  |}^2_{L^2}\\
\\
=\int_{{\mathbb R}^D}\left( {\cal F}\left( v^{\rho,r,\nu}_i(t_0,.)\right)(\xi) \exp\left(-4\pi^2\rho r^2\nu \tau |\xi|^2 \right)\right) ^2d\xi\\
\\
= \int_{{\mathbb R}^D}\left( {\cal F}\left( v^{\rho,r,\nu}_i(t_0,.)\right)^2(\xi) \exp\left(-8\pi^2\rho r^2\nu \tau |\xi|^2 \right)\right)d\xi\\
\\
=\int_{{\mathbb R}^D\setminus {\{|\xi_j|\leq \Delta},1\leq j\leq D\}}\left( {\cal F}\left( v^{\rho,r,\nu}_i(t_0,.)\right)^2(\xi) \exp\left(-8\pi^2\rho r^2\nu \tau |\xi|^2 \right)\right)d\xi\\
\\
+\int_{\{|\xi_j|\leq \Delta,~1\leq j\leq D\}}\left( {\cal F}\left( v^{\rho,r,\nu}_i(t_0,.)\right)^2(\xi) \exp\left(-8\pi^2\rho r^2\nu \tau |\xi|^2 \right)\right)d\xi\\
\\
\leq\int_{{\mathbb R}^D}\left( {\cal F}\left( v^{\rho,r,\nu}_i(t_0,.)\right)^2(\xi) \exp\left(-8\pi^2r^2\nu \tau \Delta^2 \right)\right)d\xi\\
\\
+{\Big |}\int_{\{|\xi_j|\leq \Delta,~1\leq j\leq D\}}{\Big (} {\cal F}\left( v^{\rho,r,\nu}_i(t_0,.)\right)^2(\xi)\times \\
\\
 \times \left( \exp\left(-8\pi^2\rho r^2\nu \tau |\xi|^2 \right)-\exp\left(-8\pi^2\rho r^2\nu \tau \Delta^2_0 \right)\right){\Big )} d\xi{\Big |}\\
\\
\leq {\big |} {\cal F} (v^{\rho,r,\nu}_i)(t_0,.){\big |}_{L^2}^2\exp\left(-8\pi^2\rho r^2\nu \tau \Delta^2_0 \right)
+c^{\Delta}_n
\left(8D\pi^2\rho r^2\nu \tau \Delta_0^{2+D}\right) .
\end{array}
\end{equation}
Here, we use the assumption that $\Delta_0>0$ is small enough such that the upper bound estimate is a straightforward consequence of a Taylor formula (we may choose $\Delta_0$ small and $8\pi^2\rho r^2\nu \tau \Delta_0\leq 1$) and use the abbreviation
\begin{equation}
c^{\Delta}_n:=\sup_{\{|\xi_i|\leq \Delta\}}
{\big |}{\cal F}\left( D^{\beta}_xv^{\rho,r,\nu}_i(t_0,.)\right)^2(\xi){\big |}.
\end{equation}
which is a finite constant (since ${\big |} v^{\nu}_i(t_0,.){\big |}_{H^2\cap C^2}$ is finite).

If we take the square root we may use the asymptotics $\sqrt{1+a}=1+\frac{1}{2}a+O(a^2)$.

For $\tau\in [0,\Delta_0]$ and
\begin{equation}
0<\Delta_0 \leq \max\left\lbrace \frac{1}{8\pi^2r^2\nu \max\{c^{\Delta}_n,1\}},\frac{1}{2}\right\rbrace 
\end{equation}
we get (the generous) estimate
\begin{equation}\label{vest11}
\begin{array}{ll}
{\big |}v^{\rho,r,\nu}_i(t_0,.)\ast_{sp}G_{\nu}(\tau,.){\big  |}_{L^2}\leq {\big |} {\cal F} (v^{r,\nu}_i)(t_0,.){\big |}_{L^2}\exp\left(-4\pi^2\nu\rho r^2 \tau \Delta_0^2 \right)\\
\\
+c^{\Delta}_n
\left(8D\pi^2\rho r^2\nu \tau \Delta^{1+D}\right) \\
\\
\leq {\big |}  v^{r,\nu}_i(t_0,.){\big |}_{L^2}\exp\left(-4\pi^2\nu\rho r^2 \tau \Delta^2 \right)+c^{\Delta}_n
\left(8D\pi^2\rho r^2\nu \tau \Delta_0^{1+D}\right) .
\end{array}
\end{equation}
If ${\big |}  v^{r,\nu}_i(t_0,.){\big |}_{L^2}$ becomes large or $\Delta >0$ is small enough compared to ${\big |}  v^{r,\nu}_i(t_0,.){\big |}_{L^2}$, then the second summand on right side of (\ref{vest11}) is small compared to the first summand. Note that we may replace $v^{r,\nu}_i(t_0,.)$ by multivariate spatial derivatives 
$D^{\beta}_xv^{r,\nu}_i(t_0,.)$ for $0\leq |\beta|\leq m$ such that an analogous estimate holds for spatial derivatives ${\big |} D^{\beta}_xv^{r,\nu}_i(t,.){\big |}_{L^2}$ for $0\leq |\beta |\leq m$. We shall observe below that this small damping after discrete time (under the assumption that some $|D^{\beta}_xv^{r,\nu}_i(t_0,.)|_{L^2}$ exceeds a certain level (say $1$) is strong enough in order to offset possible growth caused by the nonlinear terms.

\item[vii)] We compare the damping with the upper bound of the nonlinear term. 
\begin{rem}
For analytical purposes the damping estimate $\sim \nu \rho r^2\Delta_0^3$ (for a short time interval length $\Delta_0$) can even offset the growth of a rough upper bound $\sim \rho r\Delta_0$ (even the rougher upper bound $\sim \rho r\Delta_0$) as an upper bound for the possible growth cause by the nonlinear term as the factor $r$ is arbitrary large such that $\rho r^2$ can be chosen to be large compared to $\rho r$ (especially such that for any given $\Delta_0$ $\rho r^2\Delta_0^3$ is large compared to $\rho r\Delta_0$ or even compared to $\rho r$). However, we use the finer estimates above, where this may be also of computational interest. The offset of possible growth caused by the nonlinear term becomes effectives as the time interval $\Delta_0 >0$ is small enough and the spatial parameter $r>1$ (indicating deviation of the operator class from a strong semigroup property) is large enough.      
\end{rem}

Summing up the preceding argument recall that we may replace $\tau$ by $\tau -t_0$ on order to have for given $t_0\geq 0$ and $\tau \in \left[t_0,t_0+\Delta_0\right]$ 
 \begin{equation}\label{Navlerayscheme4r**}
\begin{array}{ll}
 {\big |}D^{\beta}_xv^{\rho,r,\nu}_i(\tau,x){\big |}
\leq {\big |} D^{\beta}_x v^{r,\nu}_i(t_0,.){\big |}_{L^2}\exp\left(-4\pi^2\nu\rho r^2 (\tau -t_0) \Delta_0^2 \right)\\
\\
+c^{\Delta}_n\left(8D\pi^2\rho r^2\nu (\tau -t_0) \Delta_0^{1+D}\right)
+\rho rL_m(4\pi \rho r^2 \nu)^{\delta} \left( (\tau-t_0)^{1-\delta}C^*+\epsilon\right).
\end{array}
\end{equation} 
 
Recall from (\ref{vest11}) and analogous estimates for spatial derivatives that we have a damping estimate
\begin{equation}\label{vest11beta}
\begin{array}{ll}
{\big |}D^{\beta}_xv^{\rho,r,\nu}_i(t_0,.)\ast_{sp}G_{\nu}(\tau,.){\big  |}_{L^2}
\\
\\
\leq {\big |} D^{\beta}_x v^{r,\nu}_i(t_0,.){\big |}_{L^2}\exp\left(-4\pi^2\nu\rho r^2 \tau \Delta_0^2 \right)+c^{\Delta}_n
\left(8D\pi^2\rho r^2\nu \tau \Delta_0^{1+D}\right).
\end{array}
\end{equation}
These estimates become effective for small $\Delta_0 >0$ if the norm of (some spatial derivative) of the initial data exceeds a certain threshold, i.e., if for some $\beta$ we have ${\big |} D^{\beta}_x v^{r,\nu}_i(t_0,.){\big |}_{L^2}\geq 1$. If this threshold is realised for some $0\leq |\beta|\leq m$, then the last upper bound term in (\ref{Navlerayscheme4r**}) or in (\ref{vest11beta}) is relatively small for small $\Delta_0\geq (\tau-t_0)$ . Note that this upper bound term contains no spatial variables, i.e., it behaves like a constant with respect to spatial norms.
Define
\begin{equation}
c^{\Delta}_D:=c^{\Delta}_n
\left(8D\pi^2\rho r^2\nu \tau \Delta_0^{1+D}\right)
\end{equation} 
For $\tau \in [t_0,t_0+\Delta_0]$ we have 
\begin{equation}\label{Navlerayscheme4r***a}
\begin{array}{ll}
 {\big |}D^{\beta}_xv^{\rho,r,\nu}_i(t_0+\Delta_0,.){\big |}_{L^2}\leq {\big |}  v^{r,\nu}_i(t_0,.){\big |}_{L^2}\exp\left(-4\pi^2\nu\rho r^2)(\tau -t_0)  \Delta_0^2 \right)+c^{\Delta}_D\\
 \\
+\rho rL_m(4\pi^2 \rho r^2 \nu)^{\delta} \left( \Delta_0^{1-\delta} C^*+\epsilon\right).
\end{array}
\end{equation}  
It is sufficient to consider the case $t_0=0$ as the same following estimates hold for $t_0\geq 0$ if $\tau$ is replaced by $\tau-t_0$. In this case the relation on (\ref{Navlerayscheme4r***a}) shows us that
\begin{equation}
{\big |}D^{\beta}_xv^{\rho,r,\nu}_i(\Delta_0,.){\big |}_{L^2}\leq {\big |}D^{\beta}_xv^{\rho,r,\nu}_i(0,.){\big |}_{L^2}
\end{equation}
if
\begin{equation}\label{ss}
\begin{array}{ll}
{\big |} D^{\beta}_x v^{\rho, r,\nu}_i(0,.){\big |}_{L^2}\left( \exp\left(-4\pi^2\nu\rho r^2 (\tau -t_0) \Delta_0^2 \right)-1\right) +c^{\Delta}_D\\
 \\
+\rho rL_m(4\pi^2 \rho r^2 \nu)^{\delta} \left( \Delta_0^{1-\delta}\right) C^* +\epsilon \leq 0.
\end{array}
\end{equation}
Note that the positive real number $\epsilon$ in (\ref{c*}) is of exponential decay on a small time interval. More precisely,
\begin{equation}
\begin{array}{ll}
\epsilon \sim \sum_{j=1}^D\int_{t_0}^{t_0+\Delta}\int_{{\mathbb R}^D\setminus B}|y_j||G^{\rho,r}_{\nu,j}(\sigma,y)|d\sigma dy\downarrow 0\mbox{ as }\Delta \downarrow 0,
\end{array}
\end{equation} 
hence $\epsilon$ is comparatively small as $4\pi\rho r^2\nu (\tau-t_0)\leq 4\pi\rho r^2\nu \Delta_0$ becomes small.  Furthermore, as $c^{\Delta}_D\downarrow 0$ as  $\Delta_0 \downarrow 0$ with $\Delta ^{D+1}$, and the damping factor  ${\big |}  v^{\rho, r,\nu}_i(0,.){\big |}_{L^2}\left(1- \exp\left(-4\pi^2\nu\rho r^2 \tau \Delta_0^2 \right)\right)$  is dominant for $\tau =\Delta_0$ if ${\big |}  v^{\rho, r,\nu}_i(0,.){\big |}_{L^2}\geq 1$ such that $c^{\Delta}_D$ is relatively small compared to the modulus of the main part of this damping term
 \begin{equation}\label{mainpart}
{\big |}  v^{\rho, r,\nu}_i(0,.){\big |}_{L^2}\left( -4\pi^2\nu\rho r^2 \Delta_0 \Delta_0^2 \right).
\end{equation}
Next we consider the conditions such that the modulus of the main damping part is larger than the last term (\ref{Navlerayscheme4r***a}) (the last term with factor $\epsilon$ is comparatively small and can be neglected).  Here we observe the exponents of the parameters $\rho,r,\nu$ and $\Delta_0$ in (\ref{mainpart}) compared to the exponents of the parmeters $\rho,r$ and $\Delta_0$ of the last term in (\ref{Navlerayscheme4r***a}). For $\tau=\Delta_0$ for the damping term in (\ref{mainpart}) we have the dependence
\begin{equation}\label{depmain}
\sim \nu \rho r^2\Delta_0^3,
\end{equation}
which has to be compared with the last term in (\ref{Navlerayscheme4r***a}), where we have the dependence
\begin{equation}\label{deplast}
(\rho )^{1+\delta}( r )^{1+2\delta}(\nu )^{\delta} \left( \Delta_0^{1-\delta}\right).
\end{equation}
The estimate of convolutions of a Lipschitz function with first order spatial derivatives of the Gaussian forced us to assume $\nu \rho r^2\geq 1$ which implies large $r>1$ in general. Hence we have to impose
\begin{equation}
\delta <\frac{1}{2}~\mbox{ such that }r^{1+2\delta}<r^2.
\end{equation}

 Choosing a small step size parameter $\rho$, say $\rho=\Delta_0^{\mu}$, we have
\begin{equation}\label{rho}
(\rho )^{1+\delta}( r )^{1+2\delta} \Delta_0^{1-\delta}
=\Delta_0^{\mu(1+\delta)+1-\delta}r^{1+2\delta}
\end{equation}
which has to be compared with the damping term, where for $\rho=\Delta_0^{\mu}$ the latter has the dependence
\begin{equation}
\rho r^2\Delta_0^3=\Delta_0^{\mu+3}r^2.
\end{equation}
For $\delta <0.5$ and small $\Delta_0$ (especially $\Delta_0<1$) we get the time step size condition
\begin{equation}
\mu (1+\delta)+1-\delta>\mu +3 \mbox{ iff }\mu >\frac{2+\delta}{\delta} 
\end{equation}
which implies $\mu >5$ for $\delta <0.5$. Hence any choice with $r$ as above, $\delta \in (0,0.5)$, and $\rho$ as above with $\mu >\frac{2+\delta}{\delta}$ implies a regular upper bound for $v^r_i,¸1\leq i\leq D$. Note that $r>1$ is related to a deviation from a strong semigroup contraction principle of the operator, and may be chosen large such that $\nu \rho r^2\geq 1$ holds for given $\nu >0$. Note that in any case the upper bound constructed is at discrete times $l\Delta_0, l\in {\mathbb N}$, where ${\mathbb N}$ denotes the set of natural numbers. However, using the local contraction result and the semigroup property we get a regular  upper bound for all time.

\end{itemize}  

\section{Appendix A): Proof of local time contraction}
Again remember the notation discussed in assumption c) above. We consider the local contraction Lemma \ref{contrlem} in more detail.
It is sufficient to prove the theorem for the functions $v^{\rho}_i(\tau,x)=v_i(t,x)$, where $v_{i,t}(t,x)=v^{\rho}_{i,\tau}(\tau,x)\frac{1}{\rho}$ for $\rho\tau=t$.  
We then have
\begin{equation}\label{a}
\begin{array}{ll}
v^{\rho,k}_{i}=v_i(t_0,.)\ast_{sp}G^{\rho}_{\nu}+\rho\sum_{j=1}^DF^i_{j}(v^{\rho,k-1},\nabla v^{\rho,k-1})\ast G^{\rho}_{\nu,j}\\
\\
=v_i(t_0,.)\ast_{sp}G^{\rho}_{\nu}+\rho\sum_{j=1}^DF^i_{j,j}(v^{\rho,k-1},\nabla v^{\rho,k-1})\ast G^{\rho}_{\nu},
\end{array}
\end{equation}
where $G^{\rho}_{\nu}=G^{\rho ,1}_{\nu}$ in the notation above, i.e., $r=1$.
For $0\leq |\gamma|\leq m-1$ and related  $1\leq |\beta|\leq m$, where $\beta_p=\gamma_p+1$ and $\beta_l=\gamma_l$ for $l\neq p$ we have
\begin{equation}
\begin{array}{ll}
D^{\beta}_xv^{\rho,k}_{i}=D^{\beta}_xv_i(t_0,.)\ast_{sp}G^{\rho}_{\nu}+\rho\sum_{j=1}^DD^{\gamma}_xF^i_{j,j}(v^{\rho,k-1},\nabla v^{\rho,k-1})\ast G^{\rho}_{\nu,p}\\
\\
=D^{\beta}_xv_i(t_0,.)\ast_{sp}G^{\rho}_{\nu}+\rho D^{\gamma}_xf^i(v^{\rho,k-1},\nabla v^{\rho,k-1})\ast G^{\rho}_{\nu,p}.
\end{array}
\end{equation}
According to our notation discussed in item c) above, we have $\sum_j F^i_{j,j}=f_i$ for all $1\leq i\leq D$, so for fixed time $\tau\in [t_0,t_0+\Delta_0]$ from (\ref{a}) we get
\begin{equation}
\begin{array}{ll}
{\big |}v^{\rho,k}_{i}(\tau,.)-v^{\rho,k-1}_{i}(\tau,.){\big |}_{H^m\cap C^m}\\
\\
=\sum_{0\leq |\beta|\leq m}{\big |}D^{\beta}_xv^{\rho,k}_{i}(\tau,.)-D^{\beta}_xv^{\rho,k-1}_{i}(\tau,.){\big |}_{L^2\cap C}\\
 \\
 \leq \rho{\Big |}\left( f_i(v^{\rho,k-1},\nabla v^{\rho,k-1})- f_i(v^{\rho,k-1},\nabla v^{\rho,k-2})\right)\ast G^{\rho}_{\nu}{\Big |}_{L^2\cap C}\\
 \\
 +\rho{\Big |}\left( f_i(v^{\rho,k-1},\nabla v^{\rho,k-2})- f_i(v^{\rho,k-2},\nabla v^{\rho,k-2})\right)\ast G^{\rho}_{\nu}{\Big |}_{L^2\cap C}\\
 \\
 +\sum_{0\leq |\gamma|\leq m-1} \rho{\Big |}\left(D^{\gamma}_x f_i(v^{\rho,k-1},\nabla v^{\rho,k-1})-D^{\gamma}_x f_i(v^{\rho,k-1},\nabla v^{\rho,k-2})\right)\ast G^{\rho}_{\nu,p}{\Big |}_{L^2\cap C}\\
 \\
 +\sum_{0\leq |\gamma|\leq m-1} \rho{\Big |}\left( D^{\gamma}_xf_i(v^{\rho,k-1},\nabla v^{\rho,k-2})-D^{\gamma}_x f_i(v^{\rho,k-2},\nabla v^{\rho,k-2})\right)\ast G^{\rho}_{\nu,p}{\Big |}_{L^2\cap C}.
\end{array}
\end{equation}
If $F^i_j$ has a local interpretation, then for all $0\leq |\beta| \leq m$ $D^{\beta}_xF^i_j$ is assumed to be Lipschitz continuous with respect to all arguments on compact domains, hence $D^{\gamma}_xf$ is also Lipschitz as a sum of Lipschitz functions. We may use an uniform Lipschitz constant $L$ which serves for all arguments and get
\begin{equation}\label{third}
\begin{array}{ll}
{\big |}v^{\rho,k}_{i}(\tau,.)-v^{\rho,k-1}_{i}(\tau,.){\big |}_{H^m\cap C^m}\\
 \\
 \leq \rho{\Big |}L{\big |} v^{\rho,k-1}_{i,j}(\tau,.)-  v^{\rho,k-2}_{i,j}(\tau,.){\big |}\ast G^{\rho}_{\nu}{\Big |}_{L^2\cap C}\\
 \\
 +\rho{\Big |}L{\big |} v^{\rho,k-1}_{i}(\tau,.)-  v^{\rho,k-2}_{i}(\tau,.){\big |}\ast G^{\rho}_{\nu}{\Big |}_{L^2\cap C}\\
 \\
 +\sum_{0\leq |\gamma|\leq m-1} \rho{\Big |}L{\big |} D^{\gamma}_xv^{\rho,k-1}_{i,j}(\tau,.)-  D^{\gamma}_xv^{\rho,k-2}_{i,j}(\tau,.){\big |}\ast G^{\rho}_{\nu,p}{\Big |}_{L^2\cap C}\\
 \\
 +\sum_{0\leq |\gamma|\leq m-1} \rho{\Big |}L{\big |} D^{\gamma}_xv^{\rho,k-1}_{i}(\tau,.)-  D^{\gamma}_xv^{\rho,k-2}_{i}(\tau,.){\big |}\ast G^{\rho}_{\nu,p}{\Big |}_{L^2\cap C}.
\end{array}
\end{equation}
For a global operator $F^i_j$ with Lipschitz continuity with respect to $|.|_{L^2\cap C}$ norm a similar relation as in (\ref{third}) holds. We have considered this elsewhere in the case of the Navier Stokes operator, and may reconsider this in Appendix B below. 
As an essential case we estimate a summand of the third term on the right side of (\ref{third}). We write
\begin{equation}
G^{\rho}_{\nu,j}=1_{B^D_1(0)}G^{\rho}_{\nu,j}+1_{({\mathbb R}^D\setminus B^D_1(0))}G^{\rho}_{\nu,j},
\end{equation}
where $1{B^D_1(0)}$ is the characteristic function of the ball of radius $1$ around $0$.

It has the upper bound
\begin{equation}
\begin{array}{ll}
\rho {\Big |}L{\big |} D^{\gamma}_x
v^{\rho,k-1}_{i,j}(\tau,.)-  D^{\gamma_x}v^{\rho,k-2}_{i,j}(\tau,.){\big |}\ast 1_{B^D_1(0)}G^{\rho}_{\nu,j}{\Big |}_{L^2\cap C}\\
\\
+\rho{\Big |}L{\big |} D^{\gamma}_xv^{\rho,k-1}_{i,j}(\tau,.)-  D^{\gamma_x}v^{\rho,k-2}_{i,j}(\tau,.){\big |}\ast 1_{({\mathbb R}^D\setminus B^D_1(0))}G^{\rho}_{\nu,j}{\Big |}_{L^2\cap C}\\
\\
\leq \rho LC_{1\gamma}{\big |} D^{\gamma}_xv^{\rho,k-1}_{i,j}(\tau,.)-  D^{\gamma_x}v^{\rho,k-2}_{i,j}(\tau,.){\big |},
\end{array}
\end{equation}
where
\begin{equation}
C_{1\gamma}:=\int_{{\mathbb R}^D}{\big |}1_{B^D_1(0)}G^{\rho}_{\nu,j}(\tau,x){\big |}dxdt+
{\Big |}{\cal F}\left( 1_{({\mathbb R}^D\setminus B^D_1(0))}G^{\rho}_{\nu,j}\right) {\Big |}_{L^2\cap C}.
\end{equation}
The convolutions with $G^{\rho}_{\nu}$ can be treated similarly. We may define
\begin{equation}
C_{0\gamma}:=\int_{{\mathbb R}^D}{\big |}1_{B^D_1(0)}G^{\rho}_{\nu}(\tau,x){\big |}dxdt+
{\Big |}{\cal F}\left( 1_{({\mathbb R}^D\setminus B^D_1(0))}G^{\rho}_{\nu}\right) {\Big |}_{L^2\cap C}.
\end{equation} 
Hence
\begin{equation}\label{third*}
\begin{array}{ll}
{\big |}v^{\rho,k}_{i}(\tau,.)-v^{\rho,k-1}_{i}(\tau,.){\big |}_{H^m\cap C^m}\\
 \\
 
 +\sum_{0\leq |\beta|\leq m} \rho L\left(C_{0\gamma}+C_{1\gamma} \right) {\Big |} D^{\beta}_xv^{\rho,k-1}_{i}(\tau,.)-  D^{\beta}_xv^{\rho,k-2}_{i}(\tau,.){\Big |}_{L^2\cap C}.
\end{array}
\end{equation}
We may then choose $\rho=\frac{1}{4L\left( \sum_{0\leq |\beta|\leq m} \left(C_{0\beta}+C_{1\beta} \right) \right) }$ and get the desired contraction.

\section{Appendix B): Additional remarks concerning local time contraction in case of the Navier Stokes equation}

We consider the local time contraction result in the special case of a global operator, i.e., in the case of the incompressible Navier Stokes equation in its Leray projection form. Here we may use that Lipschitz continuity of the Leray projection term function 
\begin{equation}
x\rightarrow \int_{{\mathbb R}^D}K_{D,i}(x-y)\sum_{l,m=1}^D g_{l,m}g_{m,l}(y)dy
\end{equation}
holds, if $g_l,~1\leq l\leq D$ is located in a suitable strong function space, i.e., if $g_l\in H^m\cap C^m$ for $m\geq 2$ for all $1\leq l\leq D$.
For coordinates $(\tau,y)$ with $y=rx$ and $\rho\tau=t$ the Leray projection form of the Navier Stokes equation is
\begin{equation}\label{nsrhor2}
\frac{\partial v^{\rho,r}_i}{\partial t}-\rho r^2\nu\Delta v^{\rho,r}_i+\rho r\sum_{j=1}^Dv^{\rho,r}_j\frac{\partial v^{\rho,r}}{\partial x_j}=\rho r\int_{{\mathbb R}^D}K_{D,i}(.-y)\sum_{l,m=1}^D v^{\rho,r}_{l,m}v^{\rho,r}_{m,l}(y)dy
\end{equation}
Remember: if we use viscosity damping in order to get global regular upper bounds, then we choose $r>1$. This parameter can be used to measure the deviation from a strong semigroup property with respect to some strong function space. A small parameter $\rho>0$ is a useful tool in order to prove local time contraction. Next we reformulate the local time contraction result in this context.
Assume that we have constructed a  global regular upper bound for the original equation up to time $t_0\geq 0$. Now data $v_i(t_0,.),~1\leq i\leq D$ with $v_i(t_0,.)\in H^m\cap C^m$ for some $m\geq 2$ are given. We define a time local iteration scheme $v^{\rho,r,k}_i,~1\leq i\leq D,~k\geq 0$ in a time interval $[t_0,t_0+\Delta]$, where $t_0\geq 0$ and
\begin{equation}
v^{\rho,r,0}_i(t,.)=v^{\rho,r}_i(t_0,.),~t\in [t_0,t_0+\Delta]~\mbox{(n.b. $v^0_i(0,.)=h_i(.)$)},
\end{equation}
and  $v^{\rho,r,k}_i$ is a solution of the linearized equation
\begin{equation}
\begin{array}{ll}
v^{\rho,r,k}_{i,\tau}-\rho r^2\nu \Delta v^{\rho,r,k}_i+ \rho r\sum_{j=1}^Dv^{\rho,r,k-1}_j\frac{\partial v^{\rho,r,k-1}}{\partial x_j}\\
\\
=\rho r\int_{{\mathbb R}^D}K_{D,i}(.-y)\sum_{l,m=1}^D v^{\rho,r,k-1}_{l,m}v^{\rho,r,k-1}_{m,l}(y)dy,
\end{array}
\end{equation}
where $1\leq i\leq D$. This is a variation of the scheme above, where the first member at $k=0$ was defined by a convolution with a Gaussian. For a contraction in a $H^m\cap C^m$-norm (with respect to spatial variables) this variation of a scheme is sufficient (the result may be easily reformulated for the former scheme above and is the adapted for even the stronger function spaces considered). For the parametrized incompressible Navier Stokes equation the local contraction result in (\ref{contrlem}) has the foollowing counterpart.
\begin{lem}\label{contrlemns}
Let $t_0\geq 0$ and assume that for some $m\geq 2$ and a finite constant $C_0>0$ we have
\begin{equation}
{\big |}v^{\rho,r}_i(t_0,.){\big |}_{H^m\cap C^m}\leq C_0.
\end{equation}
Then  there exists a $\Delta >0$ dependent only on dimension and on the constant ${\big |}v^{\rho,r}_i(t_0,.){\big |}_{H^m\cap C^m}$ such that on the time interval $[t_0,t_0+\Delta]$ the  functional increments  $\delta v^{\rho,r,k+1}_j=v^{\rho,r,k+1}_j-v^{\rho,r,k}_j,~1\leq j\leq D$ satisfy
\begin{equation}
\sup_{\tau\in [t_0,t_0+\Delta]}{\big |}\delta v^{\rho,r,k+1}_j(\tau,.){\big |}_{H^m\cap C^m}\leq \frac{1}{2} \sup_{\tau\in [t_0,t_0+\Delta]}{\big |}\delta v^{\rho,r,k}_j(\tau,.){\big |}_{H^m\cap C^m}
\end{equation}
and
\begin{equation}
\sup_{\tau\in [t_0,t_0+\Delta]}{\big |}\delta v^{\rho,r,1}_i(\tau,.){\big |}_{H^m\cap C^m}\leq \frac{1}{2}.
\end{equation}

\end{lem}

In order to prove this lemma (as a footnote to Lemma \ref{contrlem} and its proof in Appendix A) we first consider classical solution representations of the time local iteration scheme. We have
\begin{equation}\label{ans}
\begin{array}{ll}
v^{\rho,r,k}_{i}=v^{\rho,r}_i(t_0,.)\ast_{sp}G^{\rho,r}_{\nu}+\rho r\left( \sum_{j=1}^Dv^{\rho,r,k-1}_j\frac{\partial v^{\rho,r,k-1}}{\partial x_j}\right) \ast G^{\rho,r}_{\nu},\\
\\
+\rho r\left( \int_{{\mathbb R}^D}K_{D,i}(.-y)\sum_{l,m=1}^D v^{\rho,r,k-1}_{l,m}v^{\rho,r,k-1}_{m,l}(.,y)dy\right) \ast G^{\rho,r}_{\nu}.
\end{array}
\end{equation}
where we recall that $G^{\rho,r}_{\nu}$ is the fundamental solution of $G^{\rho,r}_{\nu,\tau}-\rho r^2 G^{\rho,r}_{\nu}=0$, considered on the time interval $\tau\in [t_0,t_0+\Delta]$. 
For $0\leq |\gamma|\leq m-1$ and related  $1\leq |\beta|\leq m$, where $\beta_p=\gamma_p+1$ and $\beta_l=\gamma_l$ for $l\neq p$ we have
\begin{equation}\label{bns}
\begin{array}{ll}
D^{\beta}_xv^{\rho,r,k}_{i}=D^{\beta}_xv^{\rho,r}_i(t_0,.)\ast_{sp}G^{\rho,r}_{\nu}+\rho rD^{\gamma}_x\left( \sum_{j=1}^Dv^{\rho,r,k-1}_j\frac{\partial v^{\rho,r,k-1}}{\partial x_j}\right) \ast G^{\rho,r}_{\nu,p},\\
\\
+\rho rD^{\gamma}_x\left( \int_{{\mathbb R}^D}K_{D,i}(.-y)\sum_{l,m=1}^D v^{\rho,r,k-1}_{l,m}v^{\rho,r,k-1}_{m,l}(.,y)dy\right) \ast G^{\rho,r}_{\nu,p}.
\end{array}
\end{equation}
Note that for $k=1$ we have 
\begin{equation}\label{ans2}
\begin{array}{ll}
\delta v^{\rho,r,1}_{i}=v^{\rho,r,1}_{i}-v^{\rho,r,0}_{i}=v^{\rho,r}_i(t_0,.)\ast_{sp}G^{\rho,r}_{\nu}-v^{\rho,r}_i(t_0,.)\\
\\
+\rho r\left( \sum_{j=1}^Dv^{\rho,r}_j(t_0,.)\frac{\partial v^{\rho,r}}{\partial x_j}(t_0,.)\right) \ast G^{\rho,r}_{\nu},\\
\\
+\rho r\left( \int_{{\mathbb R}^D}K_{D,i}(.-y)\sum_{l,m=1}^D v^{\rho,r}_{l,m}(t_0,y)v^{\rho,r}_{m,l}(t_0,y)dy\right) \ast G^{\rho,r}_{\nu}.
\end{array}
\end{equation}
Similarly, for $0\leq |\gamma|\leq m-1$ and related  $1\leq |\beta|\leq m$, where $\beta_p=\gamma_p+1$ and $\beta_l=\gamma_l$ for $l\neq p$ we have
\begin{equation}\label{bns2}
\begin{array}{ll}
\delta D^{\beta}_xv^{\rho,r,1}_{i}=D^{\beta}_xv^{\rho,r,1}_{i}-D^{\beta}_xv^{\rho,r,0}_{i}\\
\\
=D^{\beta}_xv^{\rho,r}_i(t_0,.)\ast_{sp}G^{\rho,r}_{\nu}-D^{\beta}_xv^{\rho,r,0}_{i}\\
\\
+\rho rD^{\gamma}_x\left( \sum_{j=1}^Dv^{\rho,r}_j(t_0,.)\frac{\partial v^{\rho,r}}{\partial x_j}(t_0,.)\right) \ast G^{\rho,r}_{\nu,p},\\
\\
+\rho rD^{\gamma}_x\left( \int_{{\mathbb R}^D}K_{D,i}(.-y)\sum_{l,m=1}^D v^{\rho,r}_{l,m}(t_0,y)v^{\rho,r}_{m,l}(t_0,y)dy\right) \ast G^{\rho,r}_{\nu,p}.
\end{array}
\end{equation}
As ${\big |}D^{\beta}_xv^{\rho,r}_{m}(t_0,.){\big |}_{L^2 \cap C}$ (recall definition of the main text) spitting intergals in  (\ref{ans2}) and (\ref{bns2}) in a local part and its complement and using Young inequalities and standard $L^2$ estimates it is straightforward to prove that given $\Delta>0$ and $r>0$ there is $\rho >0$   
such that
\begin{equation}
\sup_{\tau\in [t_0,t_0+\Delta]}{\big |}\delta v^{\rho,r,1}_{i}(\tau,.){\big |}\leq \frac{1}{2}.
\end{equation}
Next we may abbreviate
\begin{equation}
\begin{array}{ll}
D^{\gamma}_xf^{ns,\rho,r}_i\equiv \rho rD^{\gamma}_xf^{ns}_i = \rho rD^{\gamma}_x\left( \sum_{j=1}^Dv^{\rho,r,k-1}_j\frac{\partial v^{\rho,r,k-1}}{\partial x_j}\right),\\
\\
+\rho rD^{\gamma}_x\left( \int_{{\mathbb R}^D}K_{D,i}(.-y)\sum_{l,m=1}^D v^{\rho,r,k-1}_{l,m}v^{\rho,r,k-1}_{m,l}(y)dy\right),
\end{array}
\end{equation}
and get
\begin{equation}\label{esth1}
\begin{array}{ll}
{\big |}v^{\rho,r,k}_{i}(\tau,.)-v^{\rho,r,k-1}_{i}(\tau,.){\big |}_{H^m\cap C^m}\\
\\
=\sum_{0\leq |\beta|\leq m}{\big |}D^{\beta}_xv^{\rho,r,k}_{i}(\tau,.)-D^{\beta}_xv^{\rho,r,k-1}_{i}(\tau,.){\big |}_{L^2\cap C}\\
 \\
 \leq \rho r{\Big |}\left( f^{ns}_i(v^{\rho,k-1},\nabla v^{\rho,k-1})- f^{ns}_i(v^{\rho,k-2},\nabla v^{\rho,k-2})\right)\ast G^{\rho ,r}_{\nu}{\Big |}_{L^2\cap C}+\\
 \\
\sum_{0\leq |\gamma|\leq m-1} \rho r{\Big |}\left(D^{\gamma}_x f^{ns}_i(v^{\rho,k-1},\nabla v^{\rho,k-1})-D^{\gamma}_x f^{ns}_i(v^{\rho,k-2},\nabla v^{\rho,k-2})\right)\ast G^{\rho ,r}_{\nu,j}{\Big |}_{L^2\cap C}.
\end{array}
\end{equation}
We estimate the last term in the latter equation (the estimation of the other terms is similar). For each $0\leq |\gamma|\leq m-1$ we have
\begin{equation}\label{esth2}
\begin{array}{ll}
\rho r{\Big |}\left( D^{\gamma}_xf^{ns}_i(v^{\rho,k-1},\nabla v^{\rho,k-1})-D^{\gamma}_x f^{ns}_i(v^{\rho,k-2},\nabla v^{\rho,k-2})\right)\ast G^{\rho, r}_{\nu,p}{\Big |}_{L^2\cap C}\\
\\
=\rho r{\Big |}D^{\gamma}_x\left( \sum_{j=1}^D\delta v^{\rho,r,k-1}_j\frac{\partial v^{\rho,r,k-1}}{\partial x_j}\right)\ast G^{\rho, r}_{\nu,j}+D^{\gamma}_x\left( \sum_{j=1}^D v^{\rho,r,k-1}_j\frac{\partial \delta v^{\rho,r,k-1}}{\partial x_j}\right)\ast G^{\rho, r}_{\nu,j},\\
\\
+2\rho rD^{\gamma}_x\left( \int_{{\mathbb R}^D}K_{D,i}(.-y)\sum_{l,m=1}^D \delta v^{\rho,r,k-1}_{l,m}v^{\rho,r,k-1}_{m,l}(y)dy\right)\ast G^{\rho, r}_{\nu,p}{\Big |}_{L^2\cap C},
\end{array}
\end{equation}
where - according to our previous notation- $\delta v^{\rho,r,k}_i=v^{\rho,r,k}_i-v^{\rho,r,k-1}_i$.
Note that for some $m\geq 2$ we assumed
\begin{equation}
{\big |}v^{\rho,r,0}_i(\tau,.){\big |}_{H^m\cap C^m}={\big |}v^{\rho,r}_i(t_0,.){\big |}_{H^m\cap C^m}\leq C_0
\end{equation}
for some finite constant $C_0>0$. Assuming for $k\geq 1$
\begin{equation}\label{ck-1}
\sup_{\tau\in [t_0,t_0+\Delta]}{\big |}v^{\rho,r,k-1}_i(\tau,.){\big |}_{H^m\cap C^m}\leq C_{k-1}
\end{equation}
from (\ref{esth1}) and (\ref{esth2}) we get
\begin{equation}\label{esth1}
\begin{array}{ll}
\sup_{\tau\in [t_0,t_0+\Delta]}{\big |}\delta v^{\rho,r,k}_{i}(\tau,.){\big |}_{H^m\cap C^m}\\
\\
\leq \rho r c_{Dm}C_GC_KC_{k-1}\sup_{\tau\in [t_0,t_0+\Delta]}{\big |}\delta v^{\rho,r,k-1}_{i}(\tau,.){\big |}_{H^m\cap C^m}.
\end{array}
\end{equation}
Here $c_{Dm}$ is determined by the number of terms on the right side of (\ref{ans}) plus the number of terms (of derivative expansions) of (\ref{bns}) for all $0\leq |\beta|\leq m$. This constant may be determined by elementary combinatorics for numerical purposes (which is not our main interest here). Furthermore we may use local integrability of the Gaussian and first order spatial derivatives of the Gaussian on open time intervals (by standard Gaussian upper bounds, cf. main text above) and choose   
\begin{equation}
\begin{array}{ll}
C_G={\big |}1_{B^D_1(0)}G^{\rho,r}_{\nu}{\big |}_{ L^1\left((0,\Delta)\times {\mathbb R}^D \right) }+\sum_{j=1}^n{\big |}1_{B^D_1(0)}G^{\rho,r}_{\nu,j}{\big |}_{L^1\left((0,\Delta)\times {\mathbb R}^D \right)}\\
\\
+\sum_{p=1}^2{\big |}1_{{\mathbb R}^D\setminus B^D_1(0)}G^{\rho,r}_{\nu}{\big |}_{L^p\left((0,\Delta)\times {\mathbb R}^D \right)}+\sum_{j=1}^n{\big |}1_{{\mathbb R}^D\setminus B^D_1(0)}G^{\rho,r}_{\nu,j}{\big |}_{L^p\left((0,\Delta)\times {\mathbb R}^D \right)}.
\end{array}
\end{equation}
As for the constant $C_K$ we may also use local intergability of the first order derivatives of the Laplacian kernel and $L^2$-integrability outside a ball and choose
\begin{equation}
\begin{array}{ll}
C_K=1+\sum_{i=1}^D{\big |}1_{B^D_1(0)}K_{D,i}{\big |}_{ L^1\left((0,\Delta)\times {\mathbb R}^D \right) }\\
\\
+\sum_{i=1}^D{\big |}1_{{\mathbb R}^D\setminus B^D_1(0)}K_{D,i}{\big |}_{ L^2\left((0,\Delta)\times {\mathbb R}^D \right) }.
\end{array}
\end{equation}
Next note that 
\begin{equation}
\begin{array}{ll}
\sup_{\tau\in [t_0,t_0+\Delta]}{\big |}v^{\rho,r,k}_i(\tau,.){\big |}_{H^m\cap C^m}\leq {\big |}v^{\rho,r,0}_i(t_0,.){\big |}_{H^m\cap C^m}\\
\\
+\sum_{l=1}^k\sup_{\tau\in [t_0,t_0+\Delta]}{\big |}\delta v^{\rho,r,1}_i(\tau,.){\big |}_{H^m\cap C^m},
\end{array}
\end{equation} 
such that the choice 
\begin{equation}
\rho \leq \frac{1}{2r c_{Dm}C_GC_K(C_0+1)}
\end{equation}
leads to
\begin{equation}
\sup_{\tau\in [t_0,t_0+\Delta]}{\big |}v^{\rho,r,k}_i(\tau,.){\big |}_{H^m\cap C^m}\leq C_0+1,
\end{equation}
or $C_{k-1}\leq C_0+1$ for all $k\geq 1$ in the induction hypothesis in (\ref{ck-1}) above.

\section{Appendix C): Additional remarks concerning randomized models (appearance of supercritical barriers) from the perspective of Trotter product formulas}

Finally, we remark that additional ideas are required if this method is to be applied to randomized models, e.g., random averaged Navier Stokes equations or Navier Stokes equations with stochastic white noise. The reason is that the viscosity damping can only offset possible growth caused by the nonlinear terms if the functional increments of the Burgers term and the Leray projection term (considered over a small time interval $\left[t_0,t_0+\Delta \right]$ with $\Delta >0$ arbitrary)  live in strong spaces. For example, if the problem is posed on a torus (periodic boundary conditions), and regularity is measured by the order of decay of the modes of a local solution function as frequencies become large, then the decay of the functional increments of the Burgers and Leray terms is determined essentially by the decay  multiplication of infinite matrices of modes applied to the data  (with a diagonal weight for the first order derivative of the Laplace kernel in the case of the Leray projection term). If a noise term enters the equation or a random averaged model is considered then a direct application of the iteration scheme would lead to the effect that modes are mixed up and even properties which may be preserved pathwise may not be preserved 'in expectation'. Dependent on the model of course, supercritical barriers appear for strong solutions, and stronger damping would be needed in order to establish global existence of strong solutions. Let us consider this in more detail and consider the Navier stokes equation problem posed on a torus ${\mathbb T}^D_l$ of dimension $D$ and diameter $l$ (the diameter may be set to $l=1$ of course).

Given a spatial scaling parameter value $r>0$ and writing the velocity component $v^r_i=v^r_i(t,y)$ for fixed $t\geq 0$ in the analytic basis $\left\lbrace \exp\left( \frac{2\pi i\alpha x}{l}\right),~\alpha \in {\mathbb Z}^D\right\rbrace $
\begin{equation}
v^r_i(t,y):=\sum_{\alpha\in {\mathbb Z}^D}v_{i\alpha}(t)\exp{\left( \frac{2\pi i\alpha x}{l}\right) },
\end{equation}
the spatially scaled equation in (\ref{nsr1}) and (\ref{nsr2}) becomes an infinite ODE equation for the infinite time dependent vector function of velocity modes $v^r_{i\alpha},~\alpha\in {\mathbb Z}^D,~1\leq i\leq D$ of the form
\begin{equation}\label{navode200first}
\begin{array}{ll}
\frac{d v^r_{i\alpha}}{dt}=\sum_{j=1}^D\nu r^2 \left( -\frac{4\pi^2 \alpha_j^2}{l^2}\right)v^r_{i\alpha}
-r\sum_{j=1}^D\sum_{\gamma \in {\mathbb Z}^D}\frac{2\pi i \gamma_j}{l}v^r_{j(\alpha-\gamma)}v^r_{i\gamma}\\
\\
+r2\pi i\alpha_i1_{\left\lbrace \alpha\neq 0\right\rbrace}\frac{\sum_{j,k=1}^D\sum_{\gamma\in {\mathbb Z}^D}4\pi^2 \gamma_j(\alpha_k-\gamma_k)v^r_{j\gamma}v^r_{k(\alpha-\gamma)}}{\sum_{i=1}^D4\pi^2\alpha_i^2},
\end{array} 
\end{equation}
for all $1\leq i\leq D$ and where for all $\alpha\in {\mathbb Z}^D$ we have $v^r_{i\alpha}(0)=h_{i\alpha}$ along with the $\alpha$-modes $h_{i\alpha}$ of $h_i$. 
In a calculus with infinitesimals, we may describe the solution by Trotter product formulas obtained form iterative infinitesimal Euler steps. As an example consider nonstandard calculus (we do not enter technical peculiarities and describe 'internal sets' and 'extension of a set-theoretical universe' as these specifications are rather trivial or well-known at least).  For arbitrary $t_e>0$, a hyperfinite number $N$ and an infinitesimal $\delta t$ with $N\delta t=t_e$ we may define a nonstandard scheme
for time steps $$m\delta t \in \left\lbrace 0,\delta t,2\delta t,\cdots ,(N-1)\delta t, N\delta t=t_e\right\rbrace.$$ We have the infinitesimal Euler scheme
\begin{equation}\label{navode200first}
\begin{array}{ll}
v^r_{i\alpha}((m+1)\delta t)=v^r_{i\alpha}(m\delta t)+\sum_{j=1}^D\nu \left( -\frac{4\pi^2r^2 \alpha_j^2}{l^2}\right)v_{i\alpha}(m\delta t)\delta t\\
\\
-\sum_{j=1}^D\sum_{\gamma \in {\mathbb Z}^D}\frac{2\pi i r\gamma_j}{l}v^r_{j(\alpha-\gamma)}(m\delta t)v^r_{i\gamma}(m\delta t)\delta t\\
\\
+ 2\pi i\alpha_i1_{\left\lbrace \alpha\neq 0\right\rbrace}r\frac{\sum_{j,k=1}^D\sum_{\gamma\in {\mathbb Z}^D}4\pi^2 \gamma_j(\alpha_k-\gamma_k)v^r_{j\gamma}(m\delta t)v^r_{k(\alpha-\gamma)}(m\delta t)}{\sum_{i=1}^D4\pi^2\alpha_i^2}\delta t,
\end{array} 
\end{equation}
or with the abbreviation
\begin{equation}
\begin{array}{ll}
e^r_{ij\alpha\gamma}(m dt)=-r\frac{2\pi i (\alpha_j-\gamma_j)}{l}v^r_{i(\alpha-\gamma)}(m\delta t)\\
\\
+ r2\pi i\alpha_i1_{\left\lbrace \alpha\neq 0\right\rbrace}
4\pi^2 \frac{\sum_{k=1}^D\gamma_j(\alpha_k-\gamma_k)v^r_{k(\alpha-\gamma)}(m\delta t)}{\sum_{i=1}^D4\pi^2\alpha_i^2},
\end{array}
\end{equation}
we have
\begin{equation}\label{navode200second}
\begin{array}{ll}
v^r_{i\alpha}((m+1)\delta t)=v^r_{i\alpha}(m\delta t)+\sum_{j=1}^D\nu \left( -\frac{4\pi^2r^2 \alpha_j^2}{l^2}\right)v^r_{i\alpha}(m\delta t)\delta t\\
\\
+\sum_{j=1}^D\sum_{\gamma\in {\mathbb Z}^D}e^r_{ij\alpha\gamma}(m\delta t)v^r_{j\gamma}(m\delta t)\delta t.
\end{array} 
\end{equation}
In the case $l=1$ we have a solution representation
\begin{equation}\label{aa}
\begin{array}{ll}
\mathbf{v}^{r,F}(t_e)\doteq \Pi_{m=0}^{N-1}\left( \delta_{ij\alpha\beta}\exp\left(-\nu r^24\pi^2 \sum_{i=1}^D\alpha_i^2 \delta t \right)\right)\times\\
\\
\times \left( \exp\left( \left( \left( e^r_{ij\alpha\beta}\right)_{ij\alpha\beta}(m\delta t)\right)\delta t \right) \right) \mathbf{h}^F.
\end{array}
\end{equation}
Here the symbol $\doteq$ means that the identity holds up to an infinitesimal error, and we denote $\mathbf{v}^{r,F}=(v^{r,F}_1,\cdots v^{r,F}_D)^T$ with $D$ infinite vectors $v^{r,F}_i=(v^r_{i\alpha})_{\alpha \in {\mathbb Z}^D}$. The entries in $(\delta_{ij\alpha\beta})$ are Kronecker-$\delta$s which describe the unit $D{\mathbb Z}^D\times D{\mathbb Z}^D$-matrix. (We note that internal counterparts of the letter sets should be considered in order to draw global conclusions using the usual methods- again these extensions are rather trivial). The formula in (\ref{aa}) is easily verified via formulas at each time step $m$ of the form 
\begin{equation}
\left( \delta_{ij\alpha\beta}\exp\left(-\nu 4\pi^2 r^2 \sum_{i=1}^D\alpha_i^2 \delta t \right)\right)  \left( \exp\left( \left( \left( e^r_{ij\alpha\beta}(m\delta t)\right)_{ij\alpha\beta}\right)\delta t\right)  \right)\mathbf{v}^{r,F}(m\delta t)
\end{equation}
(as a representation for $\mathbf{v}^F((m+1)\delta t)$).
The use of explicit infinitesimals allows us to have an effective use of first order equality $\doteq$ for arbitrary finite time where on an infinitesimal time level we have simplifications of the formula in (\ref{aa}) in the sense that  
\begin{equation}\label{aabb}
\mathbf{v}^{r,F}(t_e)\doteq \Pi_{m=0}^{N-1}\left( \delta_{ij\alpha\beta}\left(1-\nu 4\pi^2r^2\sum_{i=1}^D\alpha_i^2 \delta t \right)\right)  \left( 1+ \left( \left( e^r_{ij\alpha\beta}\right)_{ij\alpha\beta}(m\delta t)\right)\delta t \right)  \mathbf{h}^F, 
\end{equation}
is also valid up to order $O(\delta t^2)$ (if $\delta t$ is infinitesimal).
Similarly, at each time step number $m$ an infinitesimal Eulerstep (with error $O(\delta t^2)$) is described by 
\begin{equation}
\left( \delta_{ij\alpha\beta}\exp\left(-\nu 4\pi^2 r^2 \sum_{i=1}^D\alpha_i^2 \delta t \right)\right)  \left( 1+ \left( \left( e^r_{ij\alpha\beta}(m\delta t)\right)_{ij\alpha\beta}\right)\delta t \right)\mathbf{v}^{r,F}(m\delta t).
\end{equation}
We observe that at each time step all modes get an viscosity damping except for the zero modes. However this zero modes can be controlled and are a minor problem. We considered several ways do deal with this. Here we consider an extended Trotter product formula where the zero modes are shifted to zero at each time step and reappear in the drift of an equivalent related problem. This is effectively a control function for the zero modes (we use a different notation here than we did elsewhere).

Define $\mathbf{v}^{c,r,F}=(v^{c,r,F}_1,\cdots v^{c,r,F}_D)^T$ with $D$ infinite vectors $v^{c,r,F}_i=(v^{c,r}_{i\alpha})_{\alpha \in {\mathbb Z}^D}$, where 
at all time steps $m\geq 0$ we have $v^{c,r}_{i\alpha}(m\delta t)=v^{r}_{i\alpha}(m\delta t)-c_{i\alpha}(m\delta t)$ where $c_{i\alpha}=0$ for $\alpha\neq 0$ and $c_{i0}(m\delta t)=-v^{r}_{i0}(m\delta t)$. We may plug in $v^{r}_{i\alpha}(m\delta t)=v^{c,r}_{i\alpha}(m\delta t)+c_{i\alpha}(m\delta t)=v^{r}_{i\alpha}(m\delta t)$ into the equation in (\ref{navode200first}) and get
\begin{equation}\label{navode200firscontr}
\begin{array}{ll}
\frac{d v^{c,r}_{i\alpha}}{dt}=\sum_{j=1}^D\nu r^2 \left( -\frac{4\pi^2 \alpha_j^2}{l^2}\right)v^{c,r}_{i\alpha}
-r\sum_{j=1}^D\sum_{\gamma \in {\mathbb Z}^D}\frac{2\pi i \gamma_j}{l}\left( v^{c,r}_{j(\alpha-\gamma)}+1_{\left\lbrace \alpha=\gamma\right\rbrace }c_{j(\alpha-\gamma)}\right) v^{c,r}_{i\gamma}\\
\\
+r2\pi i\alpha_i1_{\left\lbrace \alpha\neq 0\right\rbrace}\frac{\sum_{j,k=1}^D\sum_{\gamma\in {\mathbb Z}^D}4\pi^2 \gamma_j(\alpha_k-\gamma_k)v^{c,r}_{j\gamma}v^{c,r}_{k(\alpha-\gamma)}}{\sum_{i=1}^D4\pi^2\alpha_i^2}.
\end{array} 
\end{equation}
for all $1\leq i\leq D$ and where for all $\alpha\in {\mathbb Z}^D$, and where we have $v^{c,r}_{i\alpha}(0)=h_{i\alpha}$ for $\alpha\neq 0$ and $v^{c,r}_{i0}(0)=0$ for all $1\leq i\leq D$. In the case $l=1$ we have a solution representation
\begin{equation}\label{aanew}
\begin{array}{ll}
\mathbf{v}^{c,r,F}(t_e)\doteq \Pi_{m=0}^{N-1}\left( \delta_{ij\alpha\beta}\exp\left(-\nu r^24\pi^2 \sum_{i=1}^D\alpha_i^2 \delta t \right)\right)\times\\
\\
\times  \left( \exp\left( \left( \left( e^{c,r}_{ij\alpha\beta}\right)_{ij\alpha\beta}(m\delta t)\right)\delta t \right) \right) \mathbf{h}^{c,F},
\end{array}
\end{equation}
where $\mathbf{h}^{c,F}=v^{c,r,F}(0)$ are obtained form the data $\mathbf{h}^F$ by setting the zero modes to zero. All related formulas for uncontrolled schemes above have their obvious counterparts which are obtained by substitution of $e^{r}_{ij\alpha\beta}$ by $e^{c,r}_{ij\alpha\beta}$ etc.. Obviously we have
\begin{equation}
\begin{array}{ll}
e^{c,r}_{ij\alpha\gamma}(m dt)=-r\frac{2\pi i (\alpha_j-\gamma_j)}{l}\left( v^{c,r}_{i(\alpha-\gamma)}(m\delta t)+c_{k(\alpha-\gamma)}(m\delta t)\right) \\
\\
+ r2\pi i\alpha_i1_{\left\lbrace \alpha\neq 0\right\rbrace}
4\pi^2 \frac{\sum_{k=1}^D\gamma_j(\alpha_k-\gamma_k) v^{c,r}_{k(\alpha-\gamma)}(m\delta t) }{\sum_{i=1}^D4\pi^2\alpha_i^2},
\end{array}
\end{equation}

Now assume that the data $h_i$ are in a strong Sobolev space, lets say $h_i\in H^{D+2}$ in case of dimension $D$. The dual Sobolev norm
\begin{equation}\label{hmequiv}
h_i\in h^m\equiv h^m\left({\mathbb Z}^D\right)  \mbox{iff} \sum_{\alpha\in {\mathbb Z}^D}|h_{i\alpha}|^{2}({1+|\alpha|^{2m}})<\infty
\end{equation} 
tells us that there surely is a finite constant such that
\begin{equation}
\forall \alpha \in {\mathbb Z}^D:~|h_{i\alpha}|\leq \frac{C}{1+|\alpha|^{D+2}}<\infty.
\end{equation}
\begin{equation}
\begin{array}{ll}
\mathbf{v}^{c,r,F}((m+1)\delta t)=\left( \delta_{ij\alpha\beta}\exp\left(-\nu 4\pi^2 r^2 \sum_{i=1}^D\alpha_i^2 \delta t \right)\right)\times\\
\\
\times  \left( 1+ \left( \left( e^{c,r}_{ij\alpha\beta}(m\delta t)\right)_{ij\alpha\beta}\right)\delta t \right)\mathbf{v}^{c,r,F}(m\delta t).
\end{array}
\end{equation}
Recall our notation: $\mathbf{v}^{c,r,F}(m\delta t)=\left(v^{c,r}_1(m\delta t),\cdots,v^{c,r}_D(m\delta t) \right)$, where $v^{c,r}_j(m\delta t)=\left( v^{c,r}_{j\alpha}(m\delta t)\right)_{\alpha\in {\mathbb Z}^D}$

In the Euler scheme (one Euler step) the Euler terms can be estimated by elliptic integrals, such that for some constant $c$ (depending on dimension $D$) we have
\begin{equation}\label{cC}
2\pi(D+D^2)\sum_{\beta\in {\mathbb Z}^D}\frac{|\beta| C}{1+|\alpha-\beta|^{D+2}}\frac{ C}{1+|\beta|^{D+2}}\leq \frac{ cC^2}{1+|\alpha|^{D+3}},
\end{equation}
i.e. the order of polynomial decay is preserved (even more than preserved). Further transformation to another time scale $t=\rho \tau$ may be used together with a strong damping parameter $r>1$. The nonlinear terms then get the scaling $\rho r$ while the viscosity coefficient (the coefficient of the Laplcacian) becomes $\rho r^2 \nu$. For given $\nu >0$ the parameter $\rho$ may be chosen (along with large $r>1$) such that this latter term $\rho r^2 \nu$ may be large while $\rho r$ is small such that the scaled upper bound variation of (\ref{cC}) becomes
\begin{equation}\label{rrhocC}
\rho r2\pi(D+D^2)\sum_{\beta\in {\mathbb Z}^D}\frac{|\beta| C}{1+|\alpha-\beta|^{D+2}}\frac{ C}{1+|\beta|^{D+2}}\leq \frac{ \rho r cC^2}{1+|\alpha|^{D+3}},
\end{equation}
where $\rho r cC^2$ become small such that geometric series upper bounds for the sum of functional increments can be constructed. We have discussed such niceties elsewhere, and we do not need them if we want to construct an upper bound for fixed $T >0$. 
For the regularity of the data we assume that there is a constant $C>0$ such that for all $1\leq i\leq D$ and all $\alpha\in {\mathbb Z}^D$ we have
\begin{equation}\label{condh}
{\big |}h_{i\alpha}{\big |}\leq \frac{C}{1+|\alpha|^{D+2}}.
\end{equation}
Given any $T>0$ and data as in (\ref{condh}) we may choose (generously)
\begin{equation}
r>\frac{1}{\nu}\left( C+cC^2(1+T)\right) 
\end{equation}
and get a finite constant $C^*=\left( C+cC^2(1+T)\right) $ such that for all $1\leq i\leq D$ we have a global regular upper bound 
\begin{equation}
{\big |}v_i(t,.){\big |}_{H^{n+2}}\leq C^*(1+t),~t\in [0,T].
\end{equation}
This theory cannot be applied to stochastic versions of the Navier Stokes equation, however- due to the fact that elliptic integral  relations  as in (\ref{cC}) are sensitive to such extensions. As an example consider a Navier Stokes equation model with a stochastic force (white noise), i.e., an initial value problem of the form
\begin{equation}
\frac{\partial v_{i}}{\partial t}-\nu \Delta v_i+\sum_{j=1}^Dv_j\frac{\partial v_{i}}{\partial x_j}=\nabla_i p+\frac{dW}{dt},
\end{equation}
along with incompressibility $\sum_{i=1}^D v_{i,i}=0$ and data $v_i(0,.)=h_i\in H^m\cap C^m$ for some $m\geq 2$. Here the derivative of the Wiener process $W$ may be defined via the represention
\begin{equation}
W(t)=\sum_{n=1}^{\infty}N_n\int_0^tk_n(s)ds,
\end{equation}
where $(N_n)_{n=1}^{\infty}$ is a family of independent, identical, standard Gaussian random variables, and $(k_n)_{n=1}^{\infty}$ denotes a family of orthogonal functions define for nonnegative time (this orthogonal family is orthonormal in the case of standard white noise). We have
\begin{equation}
\frac{dW}{dt}=\sum_{n=1}^{\infty}N_n k_n(t).
\end{equation}
An extension of the Trotter formula scheme above (on an infinitesimal time scale) leads to an additional source term
\begin{equation}
W(m\delta t)-W((m-1)\delta t)=\sum_{n=1}^{\infty}N_n(k_n(m\delta t)-k_n((m-1)\delta t)).
\end{equation}
For such a scheme energy can be transported from lower to higher frequencies, and -depending on the family $(k_n)$ upper bound relations of the form (\ref{cC}) may not hold for the solution. Similar consideration apply to stochastic average models. 

\section{Appendix D): Extensions of the result to models with generalised diffusion term}

The Gaussian damping estimate in item vi) can be generalised to models with variable diffusion. Such generalisations are desirable as even hydrodynamic limits of the Boltzmann equation around the Maxwell state lead to a viscosity
\begin{equation}
\nu\sim D^*\left(v\otimes v-\frac{1}{3}|v|^2I \right). 
\end{equation}
Here $D^*$ is the Legendre dual of the Dirichlet form of the collision term of the Boltzmann equation linearized at the Maxwell equilibrium.
For models with uniform elliptic diffusion terms, i.e., terms of the form
\begin{equation}
\sum_{j,k=1}^Da^i_{ij}\frac{\partial^2 v_i}{\partial x_j\partial x_k},
\end{equation}
where for some constants $0<\lambda\leq \Lambda<\infty$
\begin{equation}
\forall (t,x)\in [0,\infty)\times {\mathbb R}^D:~\lambda |\xi|^2\leq \sum_{j,k=1}^Da^i_{ij}(t,x)\xi_i\xi_j\leq \Lambda |\xi|^2,
\end{equation}
there are Gaussian upper bounds of the fundamental solution. The estimates of item vi) can the be adapted to this situation straightforwardly using these Gaussian upper bounds.

For generalisations to highly degenerate diffusions we have to recall H\"{o}rmander's result. Let us  For positive  natural numbers $m,n$ a matrix-valued function 
\begin{equation}\label{vcoeff}
x\rightarrow (v_{ji})^{n,m}(x),~1\leq j\leq n,~0\leq i\leq m,
\end{equation}
 on ${\mathbb R}^n$, defines $m+1$ smooth vector fields 
\begin{equation}\label{vvec}
V_i=\sum_{j=1}^n v_{ji}(x)\frac{\partial}{\partial x_j},
\end{equation}
where $0\leq i\leq m$. H\"{o}rmander provided a sufficient condition for the fundamental solution on $[0,\infty)\times {\mathbb R}^n$ of 
\begin{equation}
	\label{hoer1}
	\left\lbrace \begin{array}{ll}
		\frac{\partial p}{\partial t}=\frac{1}{2}\sum_{i=1}^mV_i^2p+V_0p\\
		\\
		p(0,x;y)=\delta_y(x),
	\end{array}\right.
\end{equation}
where $\delta_y(x)=\delta(x-y)$ is the Dirac delta distribution with an argument shifted by the vector $y\in {\mathbb R}^n$. The condition is:
for all $x\in {\mathbb R}^n$  assume that
\begin{equation}
H_x={\mathbb R}^n,
\end{equation}
where
\begin{equation}\label{Hoergenx}
\begin{array}{ll}
H_x:=\mbox{span}{\Big\{} V_i(x), \left[V_j,V_k \right](x),
\left[ \left[V_j,V_k \right], V_l\right](x),\\
\\
\cdots |1\leq i\leq m,~0\leq j,k,l,\cdots \leq m {\Big \}},
\end{array}
\end{equation}
and where $\left[.,.\right]$ denotes the Lie bracket of vector fields as usual. Assume that
the coefficients of the vector fields are smooth (i.e. $C^{\infty}$) and bounded  with bounded derivatives, i.e. $v_{ji} \in C_{b}^{\infty}\left({\mathbb R}^n \right)$ (sometimes linear growth for the functions $v_{ji}$ themselves is allowed ( $v_{ji} \in C_{b,l}^{\infty}\left({\mathbb R}^n \right)$ in symbols). Indeed this makes no real difference. These strong regularity assumptions are due to the fact that i	 
Then a smooth density $p$ exists and for each nonnegative natural number $j$, and multiindices $\alpha,\beta$ there are increasing functions of time
\begin{equation}\label{constAB1}
A_{j,\alpha,\beta}, B_{j,\alpha,\beta}:[0,T]\rightarrow {\mathbb R},
\end{equation}
and functions
\begin{equation}\label{constmn1}
n_{j,\alpha,\beta}, 
m_{j,\alpha,\beta}:
{\mathbb N}\times {\mathbb N}^d\times {\mathbb N}^d\rightarrow {\mathbb N},
\end{equation}
such that 
\begin{equation}\label{pxest1}
\begin{array}{ll}
{\Bigg |}\frac{\partial^j}{\partial t^j} \frac{\partial^{|\alpha|}}{\partial x^{\alpha}} \frac{\partial^{|\beta|}}{\partial y^{\beta}}p(t,x,y){\Bigg |}\\
\\
\leq \frac{A_{j,\alpha,\beta}(t)(1+x)^{m_{j,\alpha,\beta}}}{t^{n_{j,\alpha,\beta}}}\exp\left(-B_{j,\alpha,\beta}(t)\frac{(x-y)^2}{t}\right).
\end{array}
\end{equation}
As the damping estimate is in discrete time at each time step of the scheme above the local damping estimate above in item vi) can be adapted to this situation. This is a first step in order to extend the construction scheme of global regular solutions models of the form
\begin{equation}
	\label{hoer1}
	\left\lbrace \begin{array}{ll}
		\frac{\partial v_i}{\partial t}=\frac{1}{2}\sum_{i=1}^mV_i^2v_i+V_0v_i+f_i(v,\nabla v)=0,~1\leq i\leq D\\
		\\
		v_i(0,.)=h_i.
	\end{array}\right.
\end{equation}
A stronger form of condition d) is required in order to get an upper bound for the functional increments
\begin{equation}
\delta v_i(t,.)=v_i(t,.)-\int_{{\mathbb R}^D}v_i(t_0,.)p(t,x;s,y)dy=\int_{t_0}^{t_0+\Delta}\int_{{\mathbb R}^D}f_i(v(s,y),\nabla v(s,y))p(t,x;s,y)dyds,
\end{equation}
because we have to compensate for the factor  
$(1+x)^{m_{j,\alpha,\beta}}$ in the H\"{o}rmander estimate of the density $p$ (which satisfies (\ref{hoer1})). Note that some order of spatial decay is lost in general as the (smoothed) initial data term $\int_{{\mathbb R}^D}v_i(t_0,.)p(t,x;s,y)dy$ shows. However this loss of spatial decay is not increased iteratively if the nonlinear term satisfies a strong submultiplicative property which compensates the factor $(1+x)^{m_{j,\alpha,\beta}}$. Clearly, the loss of polynomial decay depends on the order of regularity in terms of the constants $m_j$, but such a loss is not iteratively increased by the linear term as a Chapman Kolmogorov semi-group property holds for the density $p$ and since the functional increments preserve polynomial decay for an appropriate strong submultiplicative property.

\section{Conclusion (and comparison with results form CKN theory)}
The preceding argument shows that there is a considerable class of partial differential equations (subsuming the incompressible Navier Stokes equation) with spatially subhomogeneous nonlinear terms which have global solution branches in strong spaces for strong data. This conclusion contradicts the view of many authors that global regular solution branches of the Navier Stokes equation have to be due to the special structure of the equation (especially the Leray project term) if there are any such solutions. A difficulty of the construction of global regular solution branches for some members of some classes of equations  arises if the equations are considered in weak function spaces, i.e., for data which live in weak functions spaces, for equations with time-dependent external forces, or for equations with stochastic force terms or stochastic average equations. Especially, in the latter cases the strong upper bounds of weakly elliptic intergals or related multiplication rules of infinite mode matrices do not hold anymore for iteration schemes in general. Energy can then be transported from lower frequencies to higher frequencies and supercritical barriers may appear, i.e. models with stronger viscosity damping terms may be necessary in order to get global regular solution branches or, at least, to avoid the appearance of singularities.  It is well known that in the case of the incompressible Navier Stokes equation global regular solution branches are unique (this is a difference to the Euler equation as we argued elsewhere). If a global regular solution branch $v_i\in C^0\left([0,T],H^m\cap C^m\right),~1\leq i\leq D$ for arbitrary $T>0$ is given in the case of the Navier Stokes equation, then a Cornwall inequality implies uniqueness, i.e., if $\tilde{v}_i,~1\leq i\leq D$ is another solution of the incompressible Navier Stokes equation, then
\begin{equation}
{\big |}\tilde{v}(t)-v(t){\big |}^2_{L^2}\leq {\big |}\tilde{v}(0)-v(0){\big |}^2_{L^2}
\exp\left(C\int_0^t\left( {\big |}{\big |}v(s){\big |}^p_{L^4}+{\big |}{\big |}v(s){\big |}^2_{L^4}\right)ds  \right) 
\end{equation}
for some $p\geq 4$ which depends on the dimension (here $p=8$ in case of dimension $D=3$ is sufficient), and where $C>0$ is a constant which depends on the dimension $D$ and the viscosity. Here the function space $H^m\cap C^m$ is appropriate for Cauchy problems on the whole domain, where these function spaces reduce in the case of a torus, of course. Such a Cornwall inequality implies that a solution branch is unique in sufficiently strong function spaces. However we have no such inequality for the whole class considered here.  
The existence of global regular solution branches depends on a local time contraction theorem in strong spaces. Such local time contraction results in strong spaces can also be used in order to obtain global regularity and existence results as a consequence of the CKN result. In this case it is essential to show that a Leray-Hopf solution is left continuous with respect to time and in strong function spaces with respect to the spatial variables. If $t_0\times {\mathbb T}^D$ is a time slice where a potential singularity of the Haussdorff set of singularities is located, then we have argued elsewhere that a Leray Hopf solution $v_i,~1\leq i\leq  D$, which is smooth on $[t_0-\epsilon,t_0)$ is left continuous at $t_0$ such that
\begin{equation}\label{vit0}
v_i(t_0,.)\in H^m\cap C^m~,1\leq i\leq D,~\mbox{ for some $m\geq 2$.}
\end{equation}
This is the essential step of course as we have local time contraction and can combine the standard techniques of construction of a Leray Hopf solution with local time contraction and the singularity analysis. Having obtained (\ref{vit0}) we may obtain a local regular solution on a time interval $[t_0,t_0+\Delta]$ by local time contraction (cf. also Theorem \ref{thmlt} for an alternative local time result) (we may then start with this local solution and extend it applying the standard techniques of projection smoothing and the usual limit procedures. The resulting extended solution has a thins set of possible singularities at some time sections $t_l\times {\mathbb T}^D$ for $t_l\in I$ in an index set $I$ of Lebegues measure zero and $t_l>t_0+\Delta$. For the smallest such $t_l$ we may consider a singularity analysis -as we have shown elsewehere- and get a regular solution on a time interval $[t_0,t_l]$ with $v_i(t_l,.)\in H^m\cap C^m$ which is left continuous at $t_l$. We may then apply the local contraction result for the new data .by the semigroup property- and repeat the procedure). Here, we can use the local time contraction results in $H^m\cap C^m$ norms (spatially).
Note however, that we can also use local time regularity results such as
\begin{thm}\label{thmlt}
For some $t_0>0$ assume that $v_i(t_0,.)\in L^p~,1\leq i\leq D,~p>D$, $\sum_{i=1}^Dv_{i,i}(t_0,.)=0$. then there is $\Delta >0$ and a unique short time solution $v_i,~1\leq i\leq D$ on the time interval $[t_0,t_0+\Delta]$ such that
\begin{equation}\label{viloc}
v_i\in C\left(\left[t_0,t_0+\Delta\right],L^p\left({\mathbb T}^D \right)  \right)\cap C^{\infty}\left(\left(t_0,t_0+\Delta\right)\times {\mathbb T}^D \right)
\end{equation}
for all $1\leq i\leq D$.
\end{thm} 
Hence it seems that global regularity results can be obtained for data $L^p\left({\mathbb T}^D\right)$ and $p>D$ if the problem is posed on the torus. Note that this also implies that $H^1$ regularity implies global smoothness. This is due to the well-known embedding
\begin{equation}
H^{p,s}\subset H^{q,t}~\mbox{ for $q>p$, $s>t$, and }\frac{1}{p}-\frac{1}{q}=\frac{1}{D}\left(s-t\right) . 
\end{equation}
Here $H^{p,s}$ is the usual $L^p$ analog of $H^s$. For example, in case $D=3$, $p=2$ and $q=4$ we get $H^1\cap H^{4,0.25}\subset L^4$. Similarly for $D=3$ and $q=3+\epsilon$ for small epsilon we get an embedding in $L^{3+\epsilon}$. Hence it seems that even $H^1$ data are sufficient.
This results certainly is stronger then the results of this paper, but the CKN results cannot be transferred to the large class of equations which is considered here. This class (and variations of it with a more general class of viscosity terms) may be important also in a physical sense, as data imply that we have to go beyond the simple Navier Stokes equation in order to describe fluids appropriately. Furthermore the approach considered here is constructive and the regular upper bounds constructed can be a guide for algorithmic solution schemes. Finally we note that the existence of a (possible highly degenerate) diffusion term is essential for the preceding argument. Viscosity limit arguments ($\nu\downarrow 0$) require additional specific structure (such as vorticity forms in the case of Navier Stokes and Euler equations), and even then limits are of restricted regularity in general.

\end{document}